\documentclass{amsart}
\usepackage{amssymb}
\usepackage{amsmath}
\usepackage{amscd}
\usepackage{url}
\usepackage{setspace} 

\newtheorem{prop}{Proposition}[section]

\newtheorem{conj}[prop]{Conjecture}

\theoremstyle{definition}

\newcommand{\Aut}{{\mathrm {Aut}}}

\newcommand{\ord}{{\mathrm {ord}}}

\newcommand{\Hom}{{\mathrm {Hom}}}

\newcommand{\Ind}{{\mathrm {Ind}}}
\newcommand{\Norm}{{\mathrm {Norm}}}
\newcommand{\Fil}{{\mathrm {Fil}}}

\newcommand{\alg}{\mathrm{alg}}
\newcommand{\tr}{{\mathrm {tr}}}
\newcommand{\Tr}{{\mathrm {Tr}}}

\newcommand{\Sym}{{\mathrm {Sym}}}

\newcommand{\spin}{{\mathrm {spin}}}
\newcommand{\Spin}{{\mathrm {Spin}}}
\newcommand{\dR}{{\mathrm {dR}}}

\newcommand{\Frob}{{\mathrm {Frob}}}

\newcommand{\Gal}{\mathrm {Gal}}

\newcommand{\diag}{\mathrm{diag}}

\newcommand{\A}{{\mathbb A}}
\newcommand{\CC}{{\mathbb C}}
\newcommand{\C}{{\mathbb C}}
\newcommand{\RR}{{\mathbb R}}
\newcommand{\R}{{\mathbb R}}

\newcommand{\QQ}{{\mathbb Q}}
\newcommand{\Q}{{\mathbb Q}}
\newcommand{\ZZ}{{\mathbb Z}}
\newcommand{\Z}{{\mathbb Z}}

\newcommand{\HH}{{\mathfrak H}}

\newcommand{\MMM}{{\mathcal M}}

\newcommand{\h}{{\mathcal H}}
\newcommand{\VVV}{{\mathbb V}}

\newcommand{\n}{{\mathfrak n}}

\newcommand{\q}{{\mathfrak q}}

\newcommand{\FF}{{\mathbb F}}

\newcommand{\GL}{\mathrm {GL}}
\newcommand{\PGL}{\mathrm {PGL}}

\newcommand{\SL}{\mathrm {SL}}
\newcommand{\Sp}{\mathrm {Sp}}
\newcommand{\SO}{\mathrm {SO}}
\newcommand{\GSp}{\mathrm {GSp}}
\newcommand{\PGSp}{\mathrm {PGSp}}

\newcommand{\Qbar}{\overline{\mathbb Q}}

\newcommand{\rhobar}{\overline{\rho}}

\makeatletter
\let\@wraptoccontribs\wraptoccontribs
\makeatother

\begin{document}
\title[Eisenstein congruences for $\SO(4,3)$ and $\SO(4,4)$]{Eisenstein congruences for $\SO(4,3)$, $\SO(4,4)$, \\ spinor and triple product $L$-values}
\author{Jonas Bergstr\"om}
\author{Neil Dummigan}
\author{Thomas M\'egarban\'e}
\date{April 29th, 2016.}
\address{Matematiska institutionen\\ Stockholms universitet\\ 106 91 Stockholm\\Sweden.}
\email{jonasb@math.su.se}
\address{University of Sheffield\\ School of Mathematics and Statistics\\
Hicks Building\\ Hounsfield Road\\ Sheffield, S3 7RH\\
U.K.}
\email{n.p.dummigan@shef.ac.uk}
\address{Centre de Math\'ematiques Laurent Schwartz\\ \'Ecole Polytechnique\\ 91128 Palaiseau Cedex\\ France.}
\email{thomas.megarbane@polytechnique.edu}
\contrib[with an appendix* by]{Tomoyoshi Ibukiyama and Hidenori Katsurada}
\address{Department of Mathematics, Graduate School of Mathematics,
Osaka University,
Machikaneyama 1-1, Toyonaka, Osaka, 560-0043 Japan}
\email{ibukiyam@math.sci.osaka-u.ac.jp}
\address{Muroran Institute of Technology, 17-1 Mizumoto,
Muroran 050-8585, Japan}
\email{hidenori@mmm.muroran-it.ac.jp}
\thanks{*This work was supported by JSPS KAKENHI Grant Number
25247001.}

\begin{abstract}
We work out instances of a general conjecture on congruences between Hecke eigenvalues of induced and cuspidal automorphic representations of a reductive group, modulo divisors of certain critical $L$-values, in the case that the group is a split orthogonal group. We provide some numerical evidence in the case that the group is $\SO(4,3)$ and the $L$-function is the spinor $L$-function of a genus $2$, vector-valued, Siegel cusp form. We also consider the case that the group is $\SO(4,4)$ and the $L$-function is a triple product $L$-function.
\end{abstract}

\maketitle

\section{Introduction}
Ramanujan discovered the congruence $\tau(p)\equiv 1+p^{11}\pmod{691}$ (for all primes $p$), where $\Delta=\sum_{n=1}^{\infty}\tau(n)q^n=q\prod_{n=1}^{\infty}(1-q^n)^{24}$. We may view this as being a congruence between Hecke eigenvalues, for $T(p)$ acting on the cusp form $\Delta$ of weight $12$ for $\SL_2(\ZZ)$, and on the Eisenstein series $E_{12}$ of weight $12$. The modulus $691$ comes from a certain $L$-function evaluated at a critical point depending on the weight; specifically it divides the numerator of the rational number $\frac{\zeta(12)}{\pi^{12}}$. Conjecture 4.2 of \cite{BD} is a very wide generalisation of Ramanujan's congruence, to congruences between Hecke eigenvalues of automorphic representations of $G(\A)$, where $\A$ is the adele ring and $G/\Q$ is any connected, split reductive group. (The case of a group split over an imaginary quadratic field was dealt with in \cite{Du}.) On one side of the congruence is a cuspidal automorphic representation $\tilde{\Pi}$. On the other is one induced from a cuspidal automorphic representation $\Pi$ of the Levi subgroup $M$ of a maximal parabolic subgroup $P$.
The modulus of the congruence comes from a critical value of a certain $L$-function, associated to $\Pi$ and to the adjoint representation of the $L$-group $\hat{M}$ on the Lie algebra $\hat{\n}$ of the unipotent radical of the maximal parabolic subgroup $\hat{P}$ of $\hat{G}$. Starting from $\Pi$, we conjecture the existence of $\tilde{\Pi}$, satisfying the congruence. Ramanujan's congruence is an instance of the case $G=\GL_2, M=\GL_1\times\GL_1$. Harder's conjecture on congruences between genus-$1$ and genus-$2$ (vector-valued) Siegel modular forms is the case $G=\GSp_2$, $P$ the Siegel parabolic, $M\simeq \GL_1\times\GL_2$. In \cite{BD} we looked at these examples, and others involving $\GSp_3$ and $G_2$.

One main focus of this paper is the case $G=\SO(n+1,n)$ and $M\simeq \GL_1\times\SO(n,n-1)$. This is arguably the most direct generalisation of the congruences of Ramanujan and Harder, which themselves reappear as the cases $n=1$ and $n=2$, via the special isomorphisms $\SO(2,1)\simeq\PGL_2$ and $\SO(3,2)\simeq\PGSp_2$. In the case $n=1$ the modulus comes from Riemann's zeta function. In the case $n=2$ it comes from the $L$-function of a genus-$1$ cuspidal Hecke eigenform. In the case $n=3$ it comes from the spinor $L$-function of a genus-$2$ eigenform. This is quite satisfying, since it was only the standard $L$-function of such a form that appeared in \cite{BD} (in the case $G=\GSp_3, M\simeq \GL_1\times\GSp_2$).

The shape of the conjectural congruences is worked out in \S 2, and the special cases $n=1,2,3$ are examined further in \S 3. Actually, for each $n$, the conjecture also predicts congruences modulo divisors of Riemann zeta-values. In \S 4 we see how such congruences are implied by Ramanujan's, combined with a conjectured functorial lift from $\SO(2,1)\times\SO(n,n-1)$ to $\SO(n+1,n)$. In \S 5 we examine how the Bloch-Kato conjecture, combined with a construction tracing its roots back to Ribet's converse to Herbrand's theorem \cite{R}, leads from congruences to divisibility of $L$-values.

Calculations of Hecke eigenvalues by Faber and van der Geer (using counts of points mod $p$ on the moduli space of principally polarized abelian surfaces) provided much numerical evidence for Harder's conjecture (i.e. the case $n=2$), in each instance confirming the congruence for $p\leq 37$ \cite{FvdG, vdG}. In this paper we provide numerical examples to support the case $n=3$, which involves $\SO(4,3)$ and spinor $L$-values. In \S 6 we find some apparent congruences of the right shape, showing that they hold for $p\leq 53$. To get Hecke eigenvalues for cuspidal automorphic representations of $\SO(4,3)$ we use the compact form $\SO(7)$ instead, which allows for the computation of traces of Hecke operators using spaces of algebraic modular forms. We use the extensive data compiled by the third-named author, who has also actually proved one of the congruences for all $p$. These examples would support the conjecture if the prime moduli of these congruences appear in the numerators of certain ratios of critical spinor $L$-values.

In \S 7 we confirm these predictions by calculating sufficiently good numerical approximations to the critical values of the spinor $L$-values in question, using Dokchitser's algorithm \cite{Do} as implemented in the computer package Magma. For this, the first $150$ coefficients of the Dirichlet series were obtained from (genus $2$) Hecke eigenvalues computed by the first-named author, extending the work of Faber and van der Geer. In the denominators of the rightmost critical values, we sometimes find primes that are moduli for Harder's congruence. This can be explained via a global torsion term in the Bloch-Kato conjecture. We also find some more large primes in numerators, predicting more congruences, which are tested in \S 8.

In \cite{Du2}, the second-named author found just a scrap of numerical evidence for an Eisenstein congruence involving $U(2,2)$ and its Siegel parabolic subgroup, calculating the eigenvalues of one Hecke operator (and its square) on a single $2$-dimensional space of algebraic modular forms for $U(4)$. For this congruence, exact $L$-values were computed. The numerical support in the present paper, for Eisenstein congruences for $G=\SO(4,3), M\simeq \GL_1\times\SO(3,2)$, is more substantial, though for spinor $L$-values we had to resort to numerical approximation before guessing exact rational ratios by truncating continued fractions. As pointed out by J. Funke, this is the first evidence for Eisenstein congruences involving a group whose associated locally symmetric space does not have a complex structure. The case $G=\SO(5,4)$, $M\simeq \GL_2\times\SO(3,2)$, will be examined elsewhere. There one needs to approximate values of a degree $8$, $\GSp_2\times\GL_2$ $L$-function, which is even more difficult.

In \S 9, we consider the case $G=\SO(n,n), M\simeq \GL_2\times\SO(n-2,n-2)$, and work out what congruence the general conjecture predicts. In \S 10 we look at the special case $n=4$. Via the central isogeny from $\SO(2,2)$ to $\PGL_2\times\PGL_2$, we get a cuspidal automorphic representation of $\SO(2,2)$ from a pair of classical Hecke eigenforms, $g$ and $h$ (always level $1$ for us). For the $\GL_2$ factor of $M$ we use another $f$, and the predicted congruences involve critical values of the triple product $L$-function attached to $f,g$ and $h$. Such critical values can be computed exactly using the pullback to $\HH_1\times\HH_1\times\HH_1$ of a genus $3$ Eisenstein series to which certain holomorphic differential operators have been applied. Ibukiyama and Katsurada already computed some examples for \cite{IKPY}, and in an appendix to this paper they clarify the method, and compute more examples. To get Hecke eigenvalues for cuspidal automorphic representations of $\SO(4,4)$ we use the compact form $\SO(8)$. We present several numerical examples supporting the conjecture.

The second-named author thanks T. Berger for the reference to \cite{Wa}, G. Chenevier, C. Faber and J. Funke for helpful comments, and
J. Voight for encouraging him to look at orthogonal groups.

\section{The setup for $G=\SO(n+1,n), M\simeq \GL_1\times\SO(n,n-1)$}
Let
$$G=\SO(n+1,n)=\{g\in M_{2n+1}:\,\,^tg{J}g={J},\,\det(g)=1\},$$
where
$${J}=\begin{pmatrix}0_n & 0 & I_n\\0 & 2 & 0\\I_n & 0 & 0_n\end{pmatrix}.$$
This is a connected, reductive (even semi-simple) algebraic group, split over  $\Q$.
It has a maximal torus $T= \{\diag(t_1,\ldots,t_n,1,t_1^{-1},\ldots,t_n^{-1}) : t_1,\ldots,t_n \in \GL_1\}$ with character group $X^*(T)$ spanned by $\{e_1,\ldots,e_n\}$ where $e_i$ sends the element $\diag(t_1,\ldots,t_n,1,t_1^{-1},\ldots,t_n^{-1})$ to $t_i$ for $1\leq i\leq n$. The cocharacter group $X_*(T)$ is spanned by $\{f_1,\ldots,f_n\}$, where $f_1:t\mapsto \diag(t,1,\ldots,1,1,t^{-1},1,\ldots,1)$, etc.~and so $\langle e_i,f_j\rangle=\delta_{ij}$, where $\langle,\rangle:X^*(T)\times X_{*}(T)\rightarrow\Z$ is the natural pairing. We can order the roots so that the set of positive roots is $\Phi_G^+=\{e_i-e_j:\,i<j\}\cup \{e_i:\,1\leq i\leq n\}\cup \{e_i+e_j:\,i< j\}$, with simple positive roots $\Delta_G=\{e_1-e_2,e_2-e_3,\ldots,e_{n-1}-e_n,e_n\}$. The half-sum of the positive roots is $\rho_G=\frac{1}{2}((2n-1)e_1+(2n-3)e_2+\cdots +e_n)$. The Weyl group $W_G$ is generated by permutations of the $t_i$ and by inversions swapping $t_i$ with $t_i^{-1}$. The long element $w_0^G$ is the product of all the inversions.

If we choose the simple root $\alpha=e_1-e_2$, this determines a maximal parabolic subgroup $P=MN$, where $N$ is the unipotent radical and $M$ is the Levi subgroup, characterised by $\Delta_M=\Delta_G-\{\alpha\}$, and then $M\simeq \GL_1\times\SO(n,n-1)$. The positive roots occurring in the Lie algebra of $N$ are $\Phi_N=\Phi^+_G-\Phi^+_M=\{e_1-e_2,\ldots,e_1-e_n,e_1,e_1+e_2,\ldots,e_1+e_n\}$, i.e.~those positive roots whose expression as a sum of simple roots includes $\alpha$. The half-sum is $\rho_P=\frac{2n-1}{2}e_1$, and $\langle\rho_P,\check{\alpha}\rangle=\frac{2n-1}{2}$, where $\check{\alpha}$ is the coroot associated with $\alpha$. Let $\tilde{\alpha}:=\frac{1}{\langle\rho_P,\check{\alpha}\rangle}\rho_P=e_1$.

Let $\hat{G}$ be the Langlands dual group of $G$. (In our particular case, $\hat{G}\simeq \Sp_n$, a symplectic group of $2n$-by-$2n$ matrices. This is explained in more detail in \cite[\S 6]{Du}.) Then $\hat{G}$ has a maximal torus $\hat{T}$ with $X^*(\hat{T})\simeq X_*(T)$ and $X_*(\hat{T})\simeq X^*(T)$. Under these isomorphisms, roots of $\hat{G}$ become coroots of $G$, and coroots of $\hat{G}$ become roots of $G$, with $\check{\Delta}:=\{\check{\beta}:\beta\in\Delta_G\}$ mapping to a set of simple positive roots for $\hat{G}$. We can define a maximal parabolic subgroup $\hat{P}$ of $\hat{G}$, with Levi subgroup characterised by having set of simple positive roots $\check{\Delta}-\{\check{\alpha}\}$, hence identifiable with $\hat{M}$. Let $\hat{N}$ be the unipotent radical of $\hat{P}$, with Lie algebra $\hat{\n}$.

Let $\Pi'$ be a cuspidal, automorphic representation of $\SO(n,n-1)(\A)$, $\Pi=1\times\Pi'$, which is a unitary, cuspidal, automorphic representation of $M(\A)$. Let $\lambda=a_1e_2+\cdots a_{n-1}e_n$, with $a_1\geq a_2\geq \ldots\geq a_{n-1}\geq 0$ be the infinitesimal character of $\Pi'_{\infty}$ (or equally of $\Pi_{\infty}$, up to $W_M$). See \cite[\S 2]{BD} for further explanation. We shall assume that the $a_i$ are all distinct, with $a_{n-1}>0$, i.e.~that $\lambda$ is regular.

Let $B$ be a Borel subgroup of $M$ containing $T$, and for $s\in\CC$, $\chi\in X^*(T)$ and a valuation $v$, $|s\chi|_v(t):=|\chi(t)|_v^s$. For any prime $p$ such that the local component $\Pi_p$ (a representation of $M(\Q_p)$) is unramified (i.e.~spherical, with a non-zero $M(\ZZ_p)$-fixed vector), $\Pi_p$ is isomorphic to a (unitarily) parabolically induced representation $\Ind_{B(\QQ_p)}^{M(\QQ_p)}(|\chi_p|_p)$ for some $\chi_p=-[\log_p(\beta_1)e_2+\log_p(\beta_2)e_3+\cdots+\log_p(\beta_{n-1})e_n]\in X^*(T)\otimes_{\Z}i\R$. The $p$-adic valuation is normalised so that $|p|_p=p^{-1}$ and thus $|\chi_p(f_{j+1}(p))|_p=\beta_j$, for $1\leq j\leq n-1$. Note that we should have $|\beta_1|=\cdots =|\beta_{n-1}|=1$ if $\Pi'_p$ is tempered, which will be the case for us. This $\chi_p\in X^*(T)\otimes i\RR$ gives rise to $t(\chi_p)\in \hat{T}(\CC)\subset \hat{M}(\CC)$ such that, for any $\mu\in X_*(T)=X^*(\hat{T})$, $\mu(t(\chi_p))=|\chi_p(\mu(p))|_p$. Thus $t(\chi_p)=\diag(1,\beta_1,\ldots,\beta_{n-1},1,\beta_1^{-1},\ldots,\beta_{n-1}^{-1})$. The conjugacy class of $t(\chi_p)$ in $\hat{M}(\C)$ is the Satake parameter of $\Pi_p$, but we shall give $\chi_p$ the same title.

Given a representation $r:\hat{M}\rightarrow \GL_d$, we may define a local $L$-factor
$$L_p(s,\Pi_p,r):=\det(I-r(t(\chi_p))p^{-s})^{-1},$$
then an $L$-function (in general incomplete)

$$L_{\Sigma}(s,\Pi,r):=\prod_{p\notin \Sigma}L_p(s,\Pi_p,r),$$
where $\Sigma$ is a finite set of primes containing all those such that $\Pi_p$ is ramified. In particular, we take for $r$ the adjoint representation of $\hat{M}$ on $\hat{\n}$, which is a direct sum of subspaces on which $\hat{T}$ acts by those positive roots of $\hat{G}$ that are not roots of $\hat{M}$. These are identified with the coroots $\check{\gamma}$ of $G$, as $\gamma$ runs through $\Phi_N$. It follows that
$$L_p(s,\Pi_p,r)^{-1}=\prod_{\gamma\in\Phi_N}(1-\check{\gamma}(t(\chi_p))p^{-s})=\prod_{\gamma\in\Phi_N}(1-|\chi_p(\check{\gamma}(p))|_pp^{-s}).$$
Actually, $r$ is a direct sum of irreducible representations $r_i$ for some $1\leq i\leq m$, where $r_i$ acts on the direct sum $\hat{\n}_i$ of root spaces for $\Phi_N^i:=\{\check{\gamma}\in\Phi_N:\langle\tilde{\alpha},\check{\gamma}\rangle=i\}$, and
$$L_{\Sigma}(s,\Pi,r)=\prod_{i=1}^mL_{\Sigma}(s,\Pi,r_i).$$ In our case $m=2$, with $\Phi_N^1=\{e_1\pm e_{j+1}: 1\leq j\leq n-1\}$ and $\Phi_N^2=\{e_1\}$.
\vskip10pt
\begin{tabular}{|c|c|c|c|}\hline $\gamma\in\Phi_N$ & $\check{\gamma}$ & $\langle\lambda+s\tilde{\alpha},\check{\gamma}\rangle$ &  $|\chi_p(\check{\gamma}(p))|_p$\\\hline  $e_1-e_{j+1}$ ($1\leq j\leq n-1$) & $f_1-f_{j+1}$ & $-a_{j}+s$ & $\beta_{j}^{-1}$\\$e_1+e_{j+1}$ ($1\leq j\leq n-1$) & $f_1+f_{j+1}$ & $a_{j}+s$ & $\beta_{j}$\\$e_1$ & $2f_1$ & $2s$ & $1$\\\hline \end{tabular}
\vskip10pt
Using the table, $L_p(s,\Pi_p,r_1)=\prod_{i=1}^{n-1}[(1-\beta_ip^{-s})(1-\beta_i^{-1}p^{-s})]$, and $L_{\Sigma}(s,\Pi,r_1)$ is the $L$-function associated with $\Pi'$ and the standard $(2n-2)$-dimensional representation of $\widehat{\SO}(n,n-1)=\Sp_{n-1}$, while $L_p(s,\Pi_p,r_2)=(1-p^{-s})$, so $L_{\Sigma}(s,\Pi,r_2)=\zeta_{\Sigma}(s)$.

For $s>0$, we consider a certain parabolically induced representation $\Ind_P^G(\Pi\otimes|s\tilde{\alpha}|)$ of $G(\A)$, which has infinitesimal character (at $\infty$)
$$\lambda+s\tilde{\alpha}=se_1+a_1e_2+\cdots a_{n-1}e_n,$$
(up to $W_G$-action). We need $s\in\frac{1}{2}+\Z$ for $L_{\Sigma}(1+2s,\Pi,r_2)$ to be critical. Then we need all the $a_i$ to be in $\frac{1}{2}+\Z$ for $\lambda+s\tilde{\alpha}$ to be algebraically integral, i.e for $\langle \lambda+s\tilde{\alpha}, \check{\beta}\rangle\in\Z$ for all $\beta\in\Phi^+_G$. (This is already true for $\beta\in\Phi^+_M$, and we can check the above table for $\beta=\gamma\in\Phi_N$.) As in \cite[\S 3]{BD}, we assume that $L_{\Sigma}(s,\Pi,r_1)$ is the value at $0$ of the $L$-function attached to a motive (or at least a premotivic structure) $\MMM(r_1,\Pi\otimes |s\tilde{\alpha}|)$. Then, for the obvious choice of $w\in W_G$,
$w(\lambda+s\tilde{\alpha})=a_1e_1+\cdots a_{n-1}e_{n-1}+se_n$, which is dominant and regular if we add the condition $s<a_{n-1}$ to those already imposed. This coincides with the condition for $L_{\Sigma}(1+s,\Pi,r_1)$ to be critical. (See the end of \cite[\S 3]{BD} for more on this.) We exclude the smallest value $s=1/2$ from the conjecture below. For $1\leq i\leq m$, dividing $L_{\Sigma}(1+is,\Pi,r_i)$ by a Deligne period, we get an algebraic number, according to Deligne's conjecture on critical values of $L$-functions \cite{De}. We shall take the Deligne period normalised as in \cite[\S 4]{BD} (see also \S 6 below), and call the algebraic number $L_{\alg,\Sigma}(1+is,\Pi,r_i)$.

Let $\h=\h(G(\QQ_p),G(\ZZ_p))$ be the Hecke algebra of $\CC$-valued, compactly supported, $G(\ZZ_p)$-bi-invariant functions on $G(\QQ_p)$. If $f\in\h$ then $f$ acts on any smooth representation of $G(\Q_p)$ by $v\mapsto \int_{G(\QQ_p)}g(v)f(g)\,dg$, where $dg$ is a left- and right-invariant Haar measure, normalised so that $G(\ZZ_p)$ has volume $1$. Then $\h$ is a commutative ring under convolution of functions (which corresponds to composition of operators), and is generated by the characteristic functions $T'_{\mu}$ of double cosets $G(\ZZ_p)\mu(p) G(\ZZ_p)$, where $\mu\in X_*(T)$ is any cocharacter. If the representation is spherical, with $G(\ZZ_p)$-fixed vector $v_0$, then necessarily $T'_{\mu}(v_0)$ is also fixed, but since $v_0$ is unique up to scalar multiples, $\h$ acts on $v_0$ by a character. The value of this character on any particular element of $\h$ is a ``Hecke eigenvalue''.

Suppose that $q>2\max\langle\lambda,\check{\gamma}\rangle +1=2a_1+1$, and let $\q$ be a prime divisor of $q$ in a number field sufficiently large to accommodate all the Hecke eigenvalues and normalised $L$-values we shall consider.
The main conjecture of \cite{BD} is that if $\ord_{\q}(L_{\alg,\Sigma}(1+is,\Pi,r_i))>0$
then there exists a tempered, cuspidal, automorphic representation $\tilde{\Pi}$ of $G(\A)$, unramified outside $\Sigma$, and with $\tilde{\Pi}_{\infty}$ of infinitesimal character $w(\lambda+s\tilde{\alpha})$, such that for all $p\notin\Sigma$, and all $\mu\in X_*(T)$, the eigenvalues
of $T'_{\mu}$ on $\tilde{\Pi}_p$ and $\Ind_P^G(\Pi_p\otimes|s\tilde{\alpha}|_p)$ are congruent modulo $\q$. (Actually, we scale $T'_{\mu}$ by a certain power of $p$ to make $T_{\mu}$, see below. Also, for $i=2$ we require only $q>2+2s$. For an additional technical condition, see \cite[\S 4]{BD}.)

The standard representation of $\hat{G}\simeq \Sp_n$ has highest weight $f_1$ (identifying $X^*(\hat{T})$ with $X_{*}(T)$) and complete set of weights $\{\pm f_1,\pm f_2,\ldots, \pm f_n\}$. Given that this is a single $W_G$-orbit, i.e.~that $f_1$ is a minuscule weight, we can calculate the ``right-hand-side'' of the congruence in the following way.
The Satake parameter of $\Ind_P^G(\Pi_p\otimes|s\tilde{\alpha}|_p)$ is $\chi_p+s\tilde{\alpha}=-[\log_p(\beta_1)e_2+\log_p(\beta_2)e_3+\cdots\log_p(\beta_{n-1})e_n]+se_1$. Using this we get the following.
\vskip10pt
\begin{tabular}{|c|c|}\hline $\mu$ & $|(\chi_p+s\tilde{\alpha})(\mu(p))|_p$ \\ \hline $\pm f_1$ & $p^{\pm s}$\\$\pm f_{i+1}$ ($1\leq i\leq n-1$) & $\beta_i^{\pm 1}$ \\ \hline \end{tabular}
\vskip10pt
The trace is $p^s+p^{-s}+\sum_{i=1}^{n-1}(\beta_i+\beta_i^{-1})$. We would multiply this by $p^{\langle \rho_G,f_1\rangle}=p^{(2n-1)/2}$ to get the eigenvalue for $T'_{f_1}$, but instead we multiply by $p^{\langle w(\lambda+s\tilde{\alpha}),f_1\rangle}=p^{a_1}$, to get the eigenvalue for $T_{f_1}$:
$$T_{f_1}(\Ind_P^G(\Pi_p\otimes |s\tilde{\alpha}|_p))=p^{a_1+s}+p^{a_1-s}+\sum_{i=1}^{n-1}p^{a_1}(\beta_i+\beta_i^{-1}).$$

\section{The cases $n=1,2,3$ with $i=1$}
\subsection{$\mathbf{n=1}$}
In the special case $n=1$, $\SO(2,1)\simeq \PGL_2$. This arises from the conjugation action of $\PGL_2$ on the $3$-dimensional space of trace-$0$ matrices, preserving the quadratic form given by the determinant. If $A=\begin{pmatrix} x_2 & x_1\\x_3 & -x_2\end{pmatrix}$ is such a trace-$0$ matrix, then $-2\det A=x_1x_3+2x_2^2+x_3x_1$ is the quadratic form associated with $J$. Under this isomorphism, $\diag(t_1,t_2)\in\PGL_2$ is sent to $\diag(t_1t_2^{-1},1,t_2t_1^{-1})\in\SO(2,1)$, as one readily checks by calculating the conjugation action on $A$. Hence the characters $ae_1$ of (the maximal torus of) $\SO(2,1)$ and $a(e'_1-e'_2)$ of $\PGL_2$ correspond, where $e'_i:\diag(t_1,t_2)\mapsto t_i$. In particular, looking at the infinitesimal character of $\tilde{\Pi}_{\infty}$ when $\tilde{\Pi}$ is generated by a cuspidal Hecke eigenform $f$ of weight $k'\geq 2$ and trivial character, $\frac{k'-1}{2}(e'_1-e'_2)$ corresponds to $\frac{k'-1}{2}e_1$.

We have $M\simeq\GL_1$, and since $n-1=0$, this special case does not quite fit into the above framework, in that $\Phi_N^1$ is empty, so there are no $a_i$, no $\beta_j$, no $L(s,\Pi,r_1)$, and no upper bound on $s$. Since $\Pi$ is the trivial representation of $M(\A)$ (with $\lambda=0$), we can take $\Sigma=\emptyset$ and $L(s,\Pi,r_2)$ is still $\zeta(s)$. Letting $k'>2$ be the even integer $1+2s$, $L(1+2s,\Pi,r_2)$ becomes $\zeta(k')$.
Though we do not have an $a_1$ when $n=1$, turning to $\PGL_2$ we use the scaling factor $p^{(k'-1)/2}$ and the bound $q>k'$ (as if $a_1=\frac{k'-1}{2}$). So, for a prime $q>k'$ dividing the numerator of the Bernoulli number $B_{k'}$, we predict a cuspidal Hecke eigenform $f$ of weight $k'$ (corresponding to $\lambda+s\tilde{\alpha}=se_1$ with $s=\frac{k'-1}{2}$) and level $1$ (because $\Sigma=\emptyset$) such that
$$a_p(f)\equiv p^{k'-1}+1\pmod{q}.$$
The right-hand-side is obtained from that in the previous section by omitting all the $\beta_i$-terms and putting $a_1=s=\frac{k'-1}{2}$. This conjecture is well-known to be true; the case $k'=12, q=691$ being Ramanujan's congruence. See \cite[\S5]{BD} for the same conjecture arrived at via $G=\GL_2$. The conjecture one obtains by artificially enlarging $\Sigma$ beyond its minimum is also true \cite{DF}, as anticipated by Harder \cite{H2}.

\subsection{$\mathbf{n=2}$}
In the special case $n=2$, $\SO(3,2)\simeq \PGSp_2$. This arises from the conjugation action of $\PGSp_2$ on the $5$-dimensional space of matrices
$$A=\begin{pmatrix} x_3 & x_2 & 0 & -x_1\\x_5 & -x_3 & x_1 & 0\\0 & x_4 & x_3 & x_5\\-x_4 & 0 & x_2 & -x_3\end{pmatrix}$$
such that $AJ=J\,^tA$, preserving the quadratic form $(1/2)\Tr (A^2)=x_1x_4+x_2x_5+2x_3^2+x_4x_1+x_5x_2$, which is the one associated with $J$. Under this isomorphism, $\diag(t_1,t_2,t_0t_1^{-1},t_0t_2^{-1}) \in \PGSp_2$ is sent to $\diag(t_1t_2t_0^{-1},t_1t_2^{-1},1,t_0t_1^{-1}t_2^{-1},t_2t_1^{-1}) \in \SO(3,2)$, and the characters $ae'_1+be'_2-\frac{1}{2}(a+b)e'_0$ of $\PGSp_2$ and $\frac{a+b}{2}e_1+\frac{a-b}{2}e_2$ of $\SO(3,2)$ correspond, where $e'_i:\diag(t_1,t_2,t_0t_1^{-1},t_0t_2^{-1})\mapsto t_i$. In particular, looking at the infinitesimal character of $\tilde{\Pi}_{\infty}$ when $\tilde{\Pi}$ is generated by a Siegel modular form $F$ of weight $\Sym^j\det^k$ and trivial character, with $k\geq 3$, $(j+k-1)e'_1+(k-2)e'_2-\frac{j+2k-3}{2}e'_0$ corresponds to $\frac{j+2k-3}{2}e_1+\frac{j+1}{2}e_2$.

If $\Pi$ comes from a cuspidal Hecke eigenform $f$ of weight $k'>2$ then $\lambda=\frac{k'-1}{2}e_2$ and $w(\lambda+s\tilde{\alpha})=\frac{k'-1}{2}e_1+se_2$. Fixing $j$ and $k$ so that this is $\frac{j+2k-3}{2}e_1+\frac{j+1}{2}e_2$, the right hand side of the congruence becomes $p^{j+k-1}+p^{k-2}+a_p(f)$. The left hand side will be the Hecke eigenvalue (for the operator usually called ``$T(p)$'') for a genus-$2$ cuspidal Hecke eigenform $F$ of weight $\Sym^j\det^k$, level $1$ if $f$ is, as long as $\tilde{\Pi}_{\infty}$ is holomorphic discrete series. The $L$-value $L_{\Sigma}(1+s,\Pi,r_1)$ is $L_{\Sigma}(f,1+s+\frac{k'-1}{2})=L_{\Sigma}(f,j+k)$. We recover Harder's conjecture \cite{H1,vdG}. See \cite[\S 7]{BD} for the same conjecture arrived at via $G=\GSp_2$.
\subsection{$\mathbf{n=3}$}
Let $\Pi'$ be a cuspidal, automorphic representation of $\PGSp_2$, generated by $F$ as in the previous subsection, so the infinitesimal character of $\Pi'_{\infty}$ is $(j+k-1)e'_1+(k-2)e'_2-\frac{j+2k-3}{2}e'_0$, which is $\frac{j+2k-3}{2}e_1+\frac{j+1}{2}e_2$ as a representation of $\SO(3,2)(\A)$, or rather $\frac{j+2k-3}{2}e_2+\frac{j+1}{2}e_3$ when $M\simeq \GL_1\times\SO(3,2)$ is viewed as a Levi subgroup of $G=\SO(4,3)$.
For a prime $p$ at which $\Pi'_p$ is unramified, let $\chi_p=-[\log_p(\alpha_1)e'_1+\log_p(\alpha_2)e'_2+\log_p(\alpha_0)e'_0]$ be the Satake parameter, where $\alpha_1\alpha_2\alpha_0^2=1$. Viewing $\Pi'$ as a representation of $\SO(3,2)(\A)$, this is $$-\frac{1}{2}[(\log_p(\alpha_1)+\log_p(\alpha_2))e_1+(\log_p(\alpha_1)-\log_p(\alpha_2))e_2].$$
Since $(\alpha_0\alpha_1\alpha_2)(\alpha_0\alpha_1)=\alpha_1$, while $(\alpha_0\alpha_1\alpha_2)/(\alpha_0\alpha_1)=\alpha_2$, and again looking at $M$ inside $G$, we get
$$\chi_p=-[\log_p(\alpha_0\alpha_1\alpha_2)e_2+\log_p(\alpha_0\alpha_1)e_3],$$ i.e. $\beta_1=\alpha_0\alpha_1\alpha_2$ and $\beta_2=\alpha_0\alpha_1$. Hence
$$L_p(s,\Pi_p,r_1)^{-1}=\prod_{i=1}^{2}[(1-\beta_ip^{-s})(1-\beta_i^{-1}p^{-s})]$$
$$=(1-\alpha_0\alpha_1\alpha_2p^{-s})(1-\alpha_0p^{-s})(1-\alpha_0\alpha_1p^{-s})(1-\alpha_0\alpha_2p^{-s}),$$
and we see that $L_{\Sigma}(s,\Pi,r_1)$ is the spinor $L$-function $L_{\Sigma}(s,F,\spin)$.

The conjecture predicts a congruence modulo $\q$ if $q>j+2k-2$ and $\ord_{\q}L_{\alg,\Sigma}(1+s,F,\spin)>0$, where $s\in\frac{1}{2}+\Z$ and $0<s<\frac{j+1}{2}$, excluding $s=1/2$. The infinitesimal character of $\tilde{\Pi}_{\infty}$ is $\frac{j+2k-3}{2}e_1+\frac{j+1}{2}e_2+se_3$, and the right-hand-side of the congruence (for the Hecke eigenvalues of $T_{f_1}$) is
$$p^{((j+2k-3)/2)+s}+p^{((j+2k-3)/2)-s}+p^{(j+2k-3)/2}(\alpha_0+\alpha_0\alpha_1+\alpha_0\alpha_2+\alpha_0\alpha_1\alpha_2)$$
$$=p^{((j+2k-3)/2)+s}+p^{((j+2k-3)/2)-s}+T(p)(F),$$
where $T(p)(F)$ denotes the eigenvalue for $T(p)$ acting on $F$. The equality of $p^{(j+2k-3)/2}(\alpha_0+\alpha_0\alpha_1+\alpha_0\alpha_2+\alpha_0\alpha_1\alpha_2)=p^{a_1}(\beta_1+\beta_1^{-1}+\beta_2+\beta_2^{-1})$ with $T(p)(F)$ follows by a calculation like that at the end of \S 2.

\section{The case $i=2$}
For any $n$, we have $L_{\Sigma}(s,\Pi,r_2)=\zeta_{\Sigma}(s)$. We have already seen what happens for $n=1$, so we assume now that $n\geq 2$, for which we have considered so far only $i=1$. The group $\SO(2,1)\times\SO(n,n-1)$ is an endoscopic group of $\SO(n+1,n)$, and there is a functorial lift from $\SO(2,1)(\A)\times\SO(n,n-1)(\A)$ to $\SO(n+1,n)(\A)$ (now known by work of Arthur and others \cite{A}), coming from the obvious homomorphism of $L$-groups $\theta:\Sp_1\times\Sp_{n-1}\rightarrow\Sp_n$. As in the case $n=1$, let $s=\frac{k'-1}{2}$, and suppose that $q>k'$ with $\ord_q(\zeta_{\alg,\Sigma}(k'))>0$. Then we know there exists a cuspidal automorphic representation $\Pi''$ of $\SO(2,1)(\A)$, unramified outside $\Sigma$, satisfying a congruence as above. Recalling that $\Pi=1\times\Pi'$, where $\Pi'$ is on $\SO(n,n-1)(\A)$, we need to let $\tilde{\Pi}$ be the functorial lift of $\Pi''\times\Pi'$. To see this, let $t(\Pi''_p)\in\Sp_1(\CC)$, $t(\Pi'_p)\in\Sp_{n-1}(\CC)$ and $t(\tilde{\Pi}_p)\in\Sp_n(\CC)$ be the Satake parameters at a prime $p\notin\Sigma$. Then $t(\tilde{\Pi}_p)=\theta(t(\Pi''_p),t(\Pi'_p))$, so
$\tr(t(\tilde{\Pi}_p))=\tr(t(\Pi''_p))+\tr(t(\Pi'_p))$. Scaling by $p^{a_1}$, and bearing in mind the congruence satisfied by $\Pi''$, we see that
$$T_{f_1}(\tilde{\Pi}_p)\equiv p^{a_1+s}+p^{a_1-s}+T_{f_1}(\Pi'_p)\pmod{\q},$$
as required,
where the second $T_{f_1}$ is for $\SO(n,n-1)$. Similar reasoning using the Satake isomorphism works for any $T_{\mu}$.

Note that the automorphic representation $\tilde{\Pi}$ might not have non-zero holomorphic vectors. For example if $n=2$, $\Sigma=\emptyset$
and $\Pi'$, $\Pi''$ come from cuspidal Hecke eigenforms $f$ and $g$ of level $1$, then there is no holomorphic Yoshida lift, but the automorphic representation still exists.

\section{The Bloch-Kato conjecture}
It is convenient to introduce a ``motivic normalisation'',
$$L(s,F,\Spin):=L\left(s-\frac{j+2k-3}{2},F,\spin\right),$$
where, as before, $F$ is
a cuspidal, genus $2$, Hecke eigenform of
weight $\Sym^j\det^k$. (In all our examples, the level is $1$ so $\Sigma=\emptyset$.) We shall assume the existence of a motive $M/\QQ$ (or at least a premotivic structure comprising realisations and comparison isomorphisms, as defined in \cite[1.1.1]{DFG}) such that $L(M,s)=L(s,F,\Spin)$. Let $E$ be the field of coefficients of $M$, and let $\q\mid q$ be a prime divisor in $E$. The Hodge type of $M$ is
$$\{(0,j+2k-3),(k-2,j+k-1),(j+k-1,k-2),(j+2k-3,0)\}.$$
We assume that $q>j+2k-2$.

Let $O_{\q}$ be the ring of integers of the completion $E_{\q}$, and $O_{(\q)}$ the localisation at $\q$ of the ring of integers $O_E$ of $E$. Choose an $O_{(\q)}$-lattice $T_B$ in the Betti realisation $H_B(M)$ in such a way that $T_{\q}:=T_B\otimes O_{\q}$ is a $\Gal(\Qbar/\Q)$-invariant lattice in the $\q$-adic realisation. Then choose an $O_{(\q)}$-lattice $T_{\dR}$ in the de Rham realisation
 $H_{\dR}(M)$ in such a way that
$$\VVV(T_{\dR}\otimes O_{\q})=T_{\q}$$ as $\Gal(\Qbar_q/\QQ_q)$-representations, where $\VVV$ is the version of the Fontaine-Lafaille functor used in \cite{DFG}. Since $\VVV$ only applies to filtered $\phi$-modules, where $\phi$ is the crystalline Frobenius, $T_{\dR}$ must be $\phi$-stable. Anyway, this choice ensures that the $\q$-part of the Tamagawa factor at $q$ is trivial (by \cite[Theorem 4.1(iii)]{BK}), thus simplifying the Bloch-Kato conjecture below. The condition $q>j+2k-2$ ensures that the condition (*) in \cite[Theorem 4.1(iii)]{BK} holds.

For $s\in \frac{1}{2}+\ZZ$ with $\frac{1}{2}<s<\frac{j+1}{2}$, let $t=1+s+\frac{j+2k-3}{2}$, a critical point at which we evaluate the $L$-function.
Let $M(t)$ be the corresponding Tate twist of the motive.
Let $\Omega(t)$ be a Deligne period scaled according to the above choice, i.e. the determinant of the isomorphism
$$H_B(M(t))^+\otimes\CC\simeq (H_{\dR}(M(t))/\Fil^0)\otimes\CC,$$
calculated with respect to bases of $(2\pi i)^tT_B^{(-1)^t}$ and $T_{\dR}/\Fil^t$, so well-defined up to $O_{(\q)}^{\times}$.

The following formulation of the ($\q$-part of the) Bloch-Kato conjecture, as applied to this situation, is based on \cite[(59)]{DFG} (where $\Sigma$ was non-empty, though), using the exact sequence in their Lemma 2.1.
\begin{conj}[Bloch-Kato]\label{BK}  For $s\in \frac{1}{2}+\ZZ$ with $\frac{1}{2}<s<\frac{j+1}{2}$, and $t=1+s+\frac{j+2k-3}{2}$,
$$\ord_{\q}\left(\frac{L(M,t)}{\Omega(t)}\right)$$
$$=\ord_{\q}\left(\frac{\#H^1_{f}(\Q,T_{\q}^*(1-t)\otimes (E_{\q}/O_{\q}))}{\#H^0(\Q,T_{\q}^*(1-t)\otimes(E_{\q}/O_{\q}))\#H^0(\Q,T_{\q}(t)\otimes(E_{\q}/O_{\q}))}\right).$$
\end{conj}
Here, $T_{\q}^*=\Hom_{O_{\q}}(T_{\q},O_{\q})$, with the dual action of $\Gal(\Qbar/\Q)$, and $\#$ denotes a Fitting ideal. On the right hand side, in the numerator, is a Bloch-Kato Selmer group with local conditions (unramified at $p\neq q$, crystalline at $p=q$) for all finite primes $p$.

Let $\tilde{\Pi}$ be a cuspidal automorphic representation of $\SO(4,3)$ satisfying the congruence
$$T_{f_1}(\tilde{\Pi})\equiv p^{((j+2k-3)/2)+s}+p^{((j+2k-3)/2)-s}+T(p)(F)\pmod{\q}.$$
We seek to explain why we should expect
$$\ord_{\q}\left(\frac{L(M,t)}{\Omega(t)}\right)>0$$
as a consequence of such a congruence, by producing a non-zero element in the Bloch-Kato Selmer group $H^1_{f}(\Q,T_{\q}^*(1-t)\otimes (E_{\q}/O_{\q}))$. The construction in this special case is hopefully somewhat easier to follow than the more general argument in \cite[\S 14]{BD}.
Suppose that $\tilde{\Pi}$ has stable, tempered Arthur parameter (see \S 6 below). In this case the functorial lift of $\tilde{\Pi}$ to $\GL_6(\A)$ is cuspidal (and self-dual), and there is an associated $\q$-adic Galois representation $\tilde{\rho}:\Gal(\Qbar/\QQ)\rightarrow \GL_6(E_{\q})$ (see \cite[Remark 7.6]{Sh}). For all primes $p\neq q$, $\tilde{\rho}$ is unramified at $p$, with $\tilde{\rho}(\Frob_p^{-1})$ conjugate to $p^{(j+2k-3)/2}t(\chi_p)$, and $\tilde{\rho}$ is crystalline at $q$. We assume that, as expected, $\tilde{\rho}$ is irreducible. Note that the theorem of Calegari and Gee \cite{CG}, on the irreducibility of $\q$-adic Galois representations attached to cuspidal automorphic representations of $\GL_n$, applies only for $n\leq 5$. Let $\rho_F$ be the $4$-dimensional $\q$-adic Galois representation attached to $\Pi_F$ by Weissauer \cite[Theorem I]{W1}. Suppose that its reduction $\rhobar_F$ is irreducible. Then the congruence of Hecke eigenvalues (viewed as traces of Frobenius) implies that the composition factors of $\overline{\tilde{\rho}}$ are $\rhobar_F$, $\FF_{\q}(s-\frac{j+2k-3}{2})$ and $\FF_{\q}(-s-\frac{j+2k-3}{2})$.

The details of the following sketch are very much as in \cite[\S 8]{Br}, where the representation is $4$-dimensional rather than $6$-dimensional.
If $q\nmid B_{2s+1}$ (Bernoulli number) then we can eliminate the possibility of a non-split extension of $\FF_{\q}(s-\frac{j+2k-3}{2})$ by $\FF_{\q}(-s-\frac{j+2k-3}{2})$ inside $\overline{\tilde{\rho}}$, and so we may choose a Galois-invariant $O_{\q}$-lattice in the space of $\tilde{\rho}$ in such a way as to get an extension of $\FF_{\q}(s-\frac{j+2k-3}{2})$ by $\rhobar_F$, hence a class in $H^1(\Q,\rhobar_F(\frac{j+2k-3}{2}-s))=H^1(\Q,T_{\q}^*(1-t)\otimes\FF_q)$, since $T_{\q}^*\simeq T_{\q}(j+2k-3)$. Using the irreducibilty of $\tilde{\rho}$, in the manner of Ribet \cite{R}, one can show that, changing lattices if necessary, it is possible to ensure that we are looking at a non-split extension, hence a non-zero cohomology class. It produces a non-zero class in $H^1(\Q,T_{\q}^*(1-t)\otimes (E_{\q}/O_{\q}))$, and using the fact that $\tilde{\rho}$ is unramified at all $p\neq q$, and crystalline at $q$, it is possible to show that this class lies in $H^1_{f}(\Q,T_{\q}^*(1-t)\otimes (E_{\q}/O_{\q}))$.

\section{Some experimental congruences when $n=3$}
Let $\SO(7)$ be the special orthogonal group of the $E_7$ root lattice, the even, positive-definite lattice of discriminant $2$, unique up to isomorphism.
This is a semi-simple group over $\ZZ$, and $\SO(7)(\ZZ)\simeq W(E_7)^+$, the even subgroup of the Weyl group, of order $1451520$. For $\mu=a_1e_1+a_2e_2+a_3e_3$ (in the notation of \cite[5.2]{CR}), with $a_1,a_2,a_3\in\ZZ$ and $a_1\geq a_2\geq a_3\geq 0$, let $V_{\mu}$ be the space of the complex representation $\theta_{\mu}$ of $\SO(7)$ with highest weight $\mu$, and let $\rho:=\frac{5}{2}e_1+\frac{3}{2}e_2+\frac{1}{2}e_3$. The infinitesimal character of the representation $\theta_{\mu}$ of $\SO(7)(\RR)$ is $\mu+\rho$. Let $K$ be the open compact subgroup $\prod_p \SO(7)(\ZZ_p)$ of $\SO(7)(\A_f)$, and let
$$M(V_{\mu},K):=\{f:\SO(7)(\A_f)\rightarrow V_{\mu}:$$
$$ ~f(gk)=f(g)~\forall k\in K,~f(\gamma g)=\theta_{\mu}(\gamma)(f(g))~\forall \gamma\in \SO(7)(\QQ)\}$$
be the space of $V_{\mu}$-valued algebraic modular forms with level $K$ (i.e. ``level $1$''), where $\A_f$ is the ``finite'' part of the adele ring.
Since $\#(\SO(7)(\QQ)\backslash \SO(7)(\A_f)/K)=1$, $M(V_{\mu},K)$ can be identified with the fixed subspace $V_{\mu}^{\SO(7)(\ZZ)}$.

For each (finite) prime $p$, $\SO(7)(\QQ_p)\simeq\SO(4,3)(\QQ_p)$, and the local Hecke algebras are naturally isomorphic. The third-named author has computed the trace of $T(p):=T_{f_1}$ on $M(V_{\mu},K)$ for all $p\leq 53$, and for $a_1\leq 13$. For details of this work, including much of the numerical data, see \cite{Me}. The data at the website \cite{Me2} may be somewhat more extensive. Note that $M(V_{\mu},K)$ is isomorphic to the direct sum of $1$-dimensional $K$-fixed parts $\pi_f^K$, where $\pi=\pi_{\infty}\times\pi_f$ runs through all the automorphic representations of $\SO(7)(\A)$ such that $\pi_{\infty}\simeq V_{\mu}$, which all appear with multiplicity $1$ according to Arthur's multiplicity formula. (See \cite[Proposition 3.6]{CR}, and note that $V_{\mu}$ is self-dual.) It follows that this trace of $T(p)$ on $M(V_{\mu},K)$ is a sum of Hecke eigenvalues for $T(p)$ acting on such $\pi_f^K$.
By Arthur's endoscopic classification \cite[Theorem* 3.19]{CR} and multiplicity formula \cite[Conjecture 3.30]{CR}, we can sometimes then deduce the eigenvalue of $T(p)$ on some automorphic representation $\pi$ of $\SO(7)(\A)$ with stable, tempered Arthur parameter (in the sense of \cite[3.20]{CR}). This Arthur parameter is in this case the functorial lift of $\pi$ to $\GL_6$ via the standard representation of the $L$-group $\Sp_3(\CC)$, and by Arthur's symplectic-orthogonal alternative \cite[Theorem* 3.9]{CR}, this comes also from a discrete automorphic representation of $\SO(4,3)$, whose Hecke eigenvalue for $T(p)$ is the same. Its infinitesimal character looks the same as $\mu+\rho$, except that $e_i$ now means again what it did in \S 2 above. This automorphic representation of $\SO(4,3)$ is cuspidal, not just discrete, by \cite[Theorem 4.3]{Wa}. This allows us to test congruences of the type appearing in the conjecture in \S 3.3. Note that we have used several results which were conditional at the time of writing of \cite{CR}, but are now all proved, by work of many people, cited in the proof of \cite[Th\'eor\`eme 7.3.4]{Me}. See also the footnote below \cite[VIII, Th\'eor\`eme 1.2]{CL}, and \cite{T2}.

{\bf Example 1: $\mathbf{\Delta_{25,17,11}}$. } The Arthur parameters of the cuspidal automorphic representations of $\SO(7)$ of level $1$ and infinitesimal character $(25/2)e_1+(17/2)e_2+(11/2)e_3$ are $\Delta_{25,17,11}$ and $\Delta_{25,11}\oplus\Delta_{17}$. These are taken from \cite[Table 13]{CR}, where they were conditional on assumptions including announced results of Arthur \cite[\S 9]{A}, so the correctness of this list was double-starred in the sense of \cite{CR}; see the paragraph preceding \cite[Th\'eor\`eme 1.5**]{CR}. But it is now known unconditionally. We use the notation of Chenevier and Renard (see \cite[3.18,4.4]{CR}), so $\Delta_{25,17,11}$ denotes a representation of level $1$ and infinitesimal character $(25/2)e_1+(17/2)e_2+(11/2)e_3$, with stable, tempered Arthur parameters, while $\Delta_{25,11}\oplus\Delta_{17}$ is an endoscopic lift of cuspidal automorphic representations of $\SO(3,2)$ and $\SO(2,1)$, of level $1$ and infinitesimal characters $(25/2)e_1+(11/2)e_2$ and $(17/2)e_1$, respectively, associated with a genus $2$ form of weight $(j,k)=(10,9)$ and a genus $1$ form of weight $18$, respectively.
The $\tr(T(p))$ in the table is on a $2$-dimensional $M(V_{\mu},K)$, but subtracting off an endoscopic contribution,
$T(p)(\Delta_{25,17,11})=\tr(T(p))-[p^4T(p)(\Delta_{17})+T(p)(\Delta_{25,11})]$. The $T(p)(\Delta_{25,11})$ and $T(p)(\Delta_{25,17})$ were computed by the method of Faber and van der Geer, as described in the next section. Note that $\Delta_{25,17}$ is associated with $(j,k)=(16,6)$.

\vskip10pt
\begin{tabular}{|c|c|c|c|c|}\hline $p$ & $\tr(T(p))$ & $T(p)(\Delta_{17})$ & $T(p)(\Delta_{25,11})$ & $T(p)(\Delta_{25,17,11})$ \\\hline $2$ & $-96$ & $-528$ & $1920$ & $6432$\\\hline $3$ & $-1417608$ & $-4284$ & $-1942920$ & $872316$ \\\hline $5$ & $-1379732700$ & $-1025850$ & $-263846100$ & $-474730350$\\\hline $7$ & $-19435961616$ & $3225992$ & $-17517760400$ & $-9663808008$\\\hline $11$ & $-13089901140888$ & $-753618228$ & $-9052465894296$ & $6996289229556$\\\hline\end{tabular}

\begin{tabular}{|c|c|c|c|}\hline $p$ & $T(p)(\Delta_{25,17,11})$ & $T(p)(\Delta_{25,17})$ & $-T(p)(\Delta_{25,17,11})$\\ & & & $+[T(p)(\Delta_{25,17})+p^7+p^{18}]$\\\hline $2$ & $6432$ & $3600$ & $2^4.3.5.23.\mathbf{47}$\\\hline $3$ & $872316$ & $37800$ & $2^9.3^3.5.7.17.\mathbf{47}$\\\hline $5$ & $-474730350$ & $687689100$ & $2^{11}.3.5^2.\mathbf{47}.89.5939$\\\hline $7$ & $-9663808008$ & $10132939600$ & $2^{10}.3^4.5^2.\mathbf{47}.16708873$\\\hline $11$ & $6996289229556$ & $5673394253304$ & $2^9.3.5^2.\mathbf{47}.8699.354135787$ \\\hline\end{tabular}
\vskip10pt
This data is consistent with (and strongly suggests) a congruence
$$T(p)(\Delta_{25,17,11})  \equiv T(p)(\Delta_{25,17})+p^7+p^{18}\pmod{47},$$
which is of the shape considered in \S 3.3, with $s=11/2$. In fact, we have checked the congruence for all primes $p\leq 53$, though the table does not go so far.
Since $47>j+2k-2=26$, we should, according to the previous section, expect a factor of $47$ to appear in a certain normalised spinor $L$-value.

{\bf Example 2: $\mathbf{\Delta_{25,15,5}}$. }
Again, the $\tr(T(p))$ in the table is on a $2$-dimensional space, but subtracting off an endoscopic contribution,
$$T(p)(\Delta_{25,15,5})=\tr(T(p))-[p^5T(p)(\Delta_{15})+T(p)(\Delta_{25,5})],$$
where $\Delta_{25,5}$ corresponds to $(j,k)=(4,12)$, for which the space of level $1$ genus $2$ cusp forms is $1$-dimensional.

\vskip10pt
\begin{tabular}{|c|c|c|c|c|}\hline $p$ & $\tr(T(p))$ & $T(p)(\Delta_{15})$ & $T(p)(\Delta_{25,5})$ & $T(p)(\Delta_{25,15,5})$ \\\hline $2$ & $6816$ & $216$ & $-96$ & $0$\\\hline $3$ & $-474120$ & $-3348$ & $-527688$ & $867132$ \\\hline $5$ & $145932324$ & $52110$ & $596139180$ & $-613050606$\\\hline $7$ & $49205357040$ & $2822456$ & $-3608884496$ & $5377223544$\\\hline $11$ & $3229012641000$ & $20586852$ & $3047542095144$ & $-3134062555596$\\\hline\end{tabular}

\begin{tabular}{|c|c|c|c|}\hline $p$ & $T(p)(\Delta_{25,15,5})$ & $T(p)(\Delta_{25,15})$ & $-T(p)(\Delta_{25,15,5})$\\ & & & $+[T(p)(\Delta_{25,15})+p^{10}+p^{15}]$\\\hline $2$ & $0$ & $-3696$ & $2^4.3^2.11.\mathbf{19}$\\\hline $3$ & $867132$ & $511272$ & $2^8.3^3.\mathbf{19}.107$\\\hline $5$ & $-613050606$ & $118996620$ & $2^{10}.3^2.11.\mathbf{19}.16229$\\\hline $7$ & $5377223544$ & $-82574511536$ & $2^9.3^4.11.\mathbf{19}.353.1523$\\\hline $11$ & $-3134062555596$ & $5064306707064$ & $2^8.3^2.\mathbf{19}.95611121987$\\\hline\end{tabular}
\vskip10pt
There appears to be a congruence mod $19$ (which again we have
checked for all $p\leq 53$). This is not a large prime in the sense of the previous example, but we shall still calculate the relevant ratio of spinor $L$-values in the next section, and look out for $19$.

{\bf Example 3: $\mathbf{\Delta_{23,13,5}}$. }
This time it is easier, since there is no endoscopic contribution to subtract off, and $M(V_{\mu},K)$ is $1$-dimensional. Again, there appears to be a congruence mod $19$, and this has been checked against the data for all primes $p\leq 53$. In fact, this congruence has very recently been proved unconditionally by the third-named author,
using scalar-valued algebraic modular forms for $\mathrm{O}(25)$, in the manner of Chenevier and Lannes's proof of Harder's mod $41$ congruence using $\mathrm{O}(24)$ (referred to in \S 7, Example 4 below). He found that the modulus of the congruence is in fact $5472=2^5.3^2.19$. This work will be described in detail elsewhere.
\vskip10pt
\begin{tabular}{|c|c|c|c|}\hline $p$ & $T(p)(\Delta_{23,13,5})$ & $T(p)(\Delta_{23,13})$ & $-T(p)(\Delta_{23,13,5})$\\ & & & $+[T(p)(\Delta_{23,13})+p^9+p^{14}]$\\\hline $2$ & $0$ & $-480$ & $2^5.3^3.\mathbf{19}$\\\hline $3$ & $-304668$ & $-73080$ & $2^8.3^2.5.\mathbf{19}.23$\\\hline $5$ & $874314$ & $-140727300$ & $2^{10}.3^3.\mathbf{19}.11353$\\\hline $7$ & $452588136$ & $-2247786800$ & $2^9.3^3.5.7.\mathbf{19}.43.1709$\\\hline $11$ & $-1090903017204$ & $168545586264$ & $2^8.3^3.5^2.11.\mathbf{19}.79.133543$\\\hline $13$ & $1624277793138$ & $-6595005104660$ & $2^{10}.3^3.7.\mathbf{19}.79.13525649$\\\hline\end{tabular}
\vskip10pt
The mod $19$ congruence was actually discovered in 2014 by C. Faber, at a time when we only had the Hecke
eigenvalue for $p=2$. He, G. van der Geer and the first-named author found what appears to be a motivic structure associated
with $\Delta_{23,13,5}$ (and likewise for several other representations), and produced the putative
Hecke eigenvalues for $p\leq 17$ as traces of Frobenius, by methods similar to [FvdG, BFvdG]. They agree
with our subsequent computations.

\section{Some experimental genus $2$ spinor $L$-values}
In terms of Satake parameters,
$$L_p(s,F,\Spin)^{-1}$$ $$=1-\lambda_pp^{-s}+\frac{1}{2}(\lambda_p^2-\lambda_{p^2})p^{-2s}-\lambda_pp^{j+2k-3-3s}+p^{2j+4k-6-4s},$$

where
$$\lambda_{p^r}:=p^{r(j+2k-3)/2}(\alpha_0^r+(\alpha_0\alpha_1)^r+(\alpha_0\alpha_2)^r+(\alpha_0\alpha_1\alpha_2)^r)$$
and $F$ is a genus-$2$ Siegel eigenform. Note that $\lambda_p=T(p)(F)$, is the Hecke eigenvalue for the Hecke operator $T(p)$ acting on the eigenform $F$.

Faber and van der Geer \cite{FvdG} showed how to obtain traces of Hecke operators on spaces of cusp forms, from traces of Frobenius on the cohomology of local systems on $\mathcal A_2$, the moduli space of principally polarized abelian surfaces (which is also a Siegel modular threefold). They assumed a conjecture on the endoscopic contribution to the cohomology, which has since been proven by Petersen \cite{P} (see also Weissauer \cite{We}), building on research of many people on the automorphic representations of $\GSp_2$. Their method involves computing the zeta-functions of hyperelliptic curves of genus~$2$ (and pairs of elliptic curves) whose Jacobians make up the points of $\mathcal A_2$ that are defined over $\FF_{p^r}$. See \cite[\S\S 23,24]{vdG} for an explanation of the method. For weights $j$, $k$ such that the space of genus-$2$ cusp forms is $1$-dimensional, the trace is an eigenvalue, and in fact their computations over $\FF_{p^r}$ lead directly to $\lambda_{p^r}$ in these cases. The computations of Faber and van der Geer gave the values of $\lambda_{p^r}$ for prime powers $p^r\leq 37$ in all $1$-dimensional cases, but that is not enough to give sufficiently good approximations to the spinor $L$-values we are interested in. So the first-named author extended their computations, writing a new C-program with which he could calculate $\lambda_{p^r}$ for prime powers up to $149$, and thus the first $150$ coefficients in the Dirichlet series. The computation for the single prime $p=149$ (most of which is independent of $(j,k)$) took roughly three CPU weeks (standard desktop computer). These trace computations have already been used in the previous section, in checking congruences for $p\leq 53$. The numbering of the first three examples is as in the previous section.

{\bf Example 1: $\mathbf{(j,k)=(16,6)}$. }
Using the computer package Magma, one can define $L(s,F,\Spin)$ by the command
$$L:=\text{LSeries}(26,[0,1,-4,-3],1,V:\text{Sign}:=1);$$
Here $26=j+2k-2$, the conjectural functional equation relating $L(s)$ and $L(26-s)$. Recall that the Hodge type of the conjectural motive of which $L(s,F,\Spin)$ is the $L$-function is
$$\{(0,j+2k-3),(k-2,j+k-1),(j+k-1,k-2),(j+2k-3,0)\}.$$
Hence the product of gamma factors in the conjectural functional equation is $\Gamma_{\CC}(s)\Gamma_{\CC}(s-(k-2))=\Gamma_{\CC}(s)\Gamma_{\CC}(s-4)$, where $\Gamma_{\CC}(s):=2(2\pi)^{-s}\Gamma(s)$ \cite[5.3]{De}. If $\Gamma_{\RR}(s):=\pi^{-s/2}\Gamma(s/2)$ then $\Gamma_{\CC}(s)=\Gamma_{\RR}(s)\Gamma_{\RR}(s+1)$, so this is $\Gamma_{\RR}(s)\Gamma_{\RR}(s+1)\Gamma_{\RR}(s-4)\Gamma_{\RR}(s-3)$, which is where the vector $[0,1,-4,-3]$ comes from. The conductor is $1$, $V$ is the sequence of the first $150$ coefficients of the Dirichlet series, and the sign in the conjectured functional equation is, using \cite[5.3]{De}, $i^{(j+2k-3)-0+1}i^{(j+k-1)-(k-2)+1}=(-1)^{j+k}=(-1)^k$ (since $j$ is always even).

Magma implements the algorithm described in Dokchitser's paper \cite{Do}, which evaluates an $L$-function using a rapidly converging series which depends for its validity on the conjectured functional equation, which is simultaneously tested, via a quantity CFENew($L$) which ought to be small if the test of the functional equation is well-passed. In our case it was $0$, to $30$ decimal places, and LCfRequired($L$)=$153$, giving the number of coefficients of the Dirichlet series that would be required to ensure $30$-digit accuracy in the evaluations. Thus our $150$ coefficients should give a good approximation. Since $L(1+s,F,\spin)=L(1+s+\frac{j+2k-3}{2},F,\Spin)$, we want $L_{\alg}((25/2)+s+1,F,\Spin)$, with $s=11/2$, i.e. $L_{\alg}(19,F,\Spin)$, where $L_{\alg}(19,F,\Spin)=L(19,F,\Spin)/\Omega(19)$, for a certain Deligne period $\Omega(19)$, to which we have no direct access. This Deligne period is the determinant of a $2$ by $2$ matrix ($2$ being half the rank of the motive), whose entries are scaled by $(2\pi i)^r$ when we make a Tate twist by an even integer $r$ (sufficiently small to stay within the critical range). Hence if we look at $\frac{L(19,F,\Spin)}{\pi^4 L(17,F,\Spin)}$ it should (up to a power of $2$) be the same as $\frac{L_{\alg}(19,F,\Spin)}{L_{\alg}(17,F,\Spin)}$, which should still have the factor of $47$ we expect in the numerator of $L_{\alg}(19,F,\Spin)$, assuming we have not been unlucky enough for it to be cancelled by any $47$ in $L_{\alg}(17,F,\Spin)$.

Using Magma we found $$\frac{L(19,F,\Spin)}{\pi^4 L(17,F,\Spin)}\approx 0.0100470823379774368182814145009.$$
Using the computer package Maple we converted this to a continued fraction
$$[0,99,1,1,7,2,6,1,6,1,877118077264803576596,1,3,2,\ldots],$$
which clearly ought to be the rational number $$[0,99,1,1,7,2,6,1,6,1]=\frac{1880}{187119}=\frac{2^3. 5. 47}{3^2. 17. 1223}.$$

For each $s\in \frac{1}{2}+\ZZ$ with $\frac{1}{2}<s<\frac{j+1}{2}$, we calculated similarly the apparent rational values $\frac{L((25/2)+s+1,F,\Spin)}{\pi^4L((25/2)+s-1,F,\Spin)}$, which are in the table below. In the second row we have listed (the Arthur parameters of) the cuspidal automorphic representations of $\SO(7)$ of level $1$ and infinitesimal character $(25/2)e_1+(17/2)e_2+se_3$. These are taken from \cite[Table 13]{CR}. Again we use the notation of Chenevier and Renard (see \cite[3.18,4.4]{CR}), so for example $\Delta_{25,17,3}^2$ denotes a pair of representations of level $1$ and infinitesimal character $(25/2)e_1+(17/2)e_2+(3/2)e_3$, with stable, tempered Arthur parameters.
\vskip10pt
\begin{tabular}{|c|c|c|c|c|}\hline $s$ & $3/2$ & $5/2$ & $7/2$ & $9/2$ \\\hline Reps. & $\Delta_{25,17,3}^2$ & $\Delta_{25,5}\oplus\Delta_{17}$ & $\Delta_{25,7}\oplus\Delta_{17}, \Delta_{25,17,7}^2$ & $\Delta_{25,9}^2\oplus\Delta_{17}$ \\\hline $\frac{L((25/2)+s+1,\Spin)}{\pi^4L((25/2)+s-1,\Spin)}$ & $\frac{5.59}{2^23^37^213}$ & $\frac{2^2}{3.5^2.7}$ & $\frac{1223}{2^23^25^359}$ & $\frac{1}{2.3.17}$ \\\hline \end{tabular}

\begin{tabular}{|c|c|c|c|}\hline $s$ & $11/2$ & $13/2$ & $15/2$\\\hline Reps. & $\Delta_{25,11}\oplus\Delta_{17},\Delta_{25,17,11}$ & $\Delta_{25,13}^2\oplus\Delta_{17}$ & $\Delta_{25,15}\oplus\Delta_{17}$\\\hline $\frac{L((25/2)+s+1,\Spin)}{\pi^4L((25/2)+s-1,\Spin)}$ & $\frac{2^3.5.47}{3^2.17.1223}$ & $\frac{2^4}{3^25^27}$ & $\frac{2^2.3}{5^2.47}$\\\hline \end{tabular}
\vskip10pt
It is striking that one sees a large prime in the numerator precisely when there is a representation $\Delta_{25,17,2s}$, with stable, tempered Arthur parameter, available to participate in the predicted congruence. These representations are the ones that should have {\em irreducible} $6$-dimensional Galois representations attached to them, leading to an explanation, via the Bloch-Kato conjecture, of the occurrence of the large prime in the $L$-value, as a consequence of the congruence, as in \S 5 above. Only in the case $s=11/2$ is there a single representation with stable, tempered Arthur parameter, so that we may easily deduce from the trace of a Hecke operator its eigenvalue for that representation, and test the predicted congruence (already done in \S 6).

{\bf Example 2: $\mathbf{(j,k)=(14,7)}$. }
\vskip10pt
\begin{tabular}{|c|c|c|c|c|c|}\hline $s$ & $5/2$ & $7/2$ & $9/2$ \\\hline Reps. & $\Delta_{25,5}\oplus\Delta_{15}, \Delta_{25,15,5}$ & $\Delta_{25,7}\oplus\Delta_{15}$ &$\Delta_{25,9}^2\oplus\Delta_{15}, \Delta_{25,15,9}$ \\\hline $\frac{L((25/2)+s+1,\Spin)}{\pi^4L((25/2)+s-1,\Spin)}$ & $\frac{2^2.19}{3.5^2.7.11}$ & $\frac{1}{2.3^2.5}$ & $\frac{557}{2.3^4.17.19}$\\\hline\end{tabular}

\begin{tabular}{|c|c|c|}\hline $s$ & $11/2$ & $13/2$\\\hline Reps. & $\Delta_{25,11}\oplus\Delta_{15}$ & $\Delta_{25,13}^2\oplus\Delta_{15}$\\\hline $\frac{L((25/2)+s+1,\Spin)}{\pi^4L((25/2)+s-1,\Spin)}$ & $\frac{2^3}{3^2.5.17}$ & $\frac{2^4.5.7}{97.557}$\\\hline\end{tabular}
\vskip10pt
We see the anticipated factor of $19$ for $s=5/2$, and the large prime $557$ for $s=9/2$ coinciding with the appearance of a stable, tempered representation that could participate in a congruence.
Since the $L$-function now vanishes at the central point, we omitted $s=3/2$, to avoid dividing by $0$. But the entry in the second row for $s=3/2$ would have been ``none'', and there is no large prime in the denominator of the first ratio of $L$-values in the third row.

We notice a large prime $97$ in the last denominator. To explain it via the Bloch-Kato conjecture (\ref{BK}), we proceed as follows. If $\rho_f$ is the $q$-adic Galois representation attached to the normalised cusp form $f$ of weight $26$ for $\SL_2(\ZZ)$ (with $\q=97$), then Harder's conjectured congruence
$$T(p)(F)\equiv a_p(f)+p^{k-2}+p^{j+k-1}\pmod{97}$$
implies that the composition factors of $\rhobar_F$ are $\FF_{97}(2-k)$, $\FF_{97}(1-j-k)$ and $\rhobar_f$. (The latter is irreducible in this case.) Note that $97$ is a divisor of $L_{\alg}(f,j+k)$.

We choose the Galois-invariant $O_{\q}$ lattice $T_{\q}$ in the space of $\rho_F$ in such a way that $\FF_{97}(1-j-k)$ is a submodule of $\rhobar_F$. Then $\FF_{97}$ is a submodule of $\rhobar_F(j+k-1)=\rhobar_F(t)$, where $t=\frac{j+2k-3}{2}+\frac{j-1}{2}+1=j+k-1$ is the rightmost critical point. (It is no accident that the exponent in the power of $p$ marks the boundary of the critical range, since it is a Hodge weight.) This contributes a factor of $97$ to the term $\#H^0(\Q,T_{\q}(t)\otimes(E_{\q}/O_{\q}))$, which (assuming it does not appear also in $\#H^1_{f}(\Q,T_{\q}^*(1-t)\otimes (E_{\q}/O_{\q}))$) should therefore appear in the denominator of the ratio $\frac{L(j+k-1,F,\Spin)}{\pi^4 L(j+k-3,F,\Spin)}$.

{\bf Example 3: $\mathbf{(j,k)=(12,7)}$. }
\vskip10pt
\begin{tabular}{|c|c|c|c|c|}\hline $s$ & $5/2$ & $7/2$ & $9/2$ & $11/2$\\\hline Reps. & $\Delta_{23,13,5}$ & none & none & none \\\hline $\frac{L((25/2)+s+1,\Spin)}{\pi^4L((25/2)+s-1,\Spin)}$ & $\frac{19}{3.5.7.13}$ & $\frac{1}{2.3^2.5}$ & $\frac{1}{5.19}$ & $\frac{2.3^2.5}{7.17.73}$\\\hline\end{tabular}
\vskip10pt
We see the anticipated factor of $19$ for $s=5/2$, and no other occurrences of representations with stable, tempered Arthur parameter, or large primes in numerators. Again, the entry in the second row for $s=3/2$ would have been ``none''. The factor of $73$ in the last denominator can be explained similarly to above.

{\bf Example 4: $\mathbf{(j,k)=(4,10)}$. }
The original numerical example of a congruence for Harder's conjecture, appearing in \cite{H1}, is
$$T(p)(\Delta_{21,5})\equiv a_p(\Delta_{21})+p^{8}+p^{13}\pmod{41}.$$
This instance of Harder's conjecture has actually been proved by Chenevier and Lannes \cite[X, Th\'eor\`eme* 4.4(1)]{CL}. We are grateful to G. Chenevier for explaining that this theorem is now unconditional, thanks to recent work of Moeglin, Waldspurger, Shelstad and Mezo. (See the paragraph following \cite[VIII, Th\'eor\`eme* 1.1]{CL}.)
We expect to see $41$ in the denominator of $\frac{L(13,F,\Spin)}{\pi^4L(11,F,\Spin)}$, which would be the only entry in the table, for $s=3/2$ ($(21/2)+s+1=13$). We find that it appears indeed to be $\frac{2}{3.41}$.

{\bf Example 5: $\mathbf{(j,k)=(18,5)}$. }
This was not suggested by any of the congruences found in the previous section, but we look at it anyway.
\vskip10pt
\begin{tabular}{|c|c|c|c|}\hline $s$ & $5/2$ & $7/2$ & $9/2$ \\\hline Reps. & $\Delta_{25,5}\oplus\Delta_{19},\Delta_{25,19,5}^2$ & $\Delta_{25,7}\oplus\Delta_{19}$ & $\Delta_{25,9}^2\oplus\Delta_{19}, \Delta_{25,19,9}^2$ \\\hline $\frac{L((25/2)+s+1,\Spin)}{\pi^4L((25/2)+s-1,\Spin)}$ & $\frac{103}{3.5.7^2.11}$ & $\frac{1}{2.3^2.5}$ & $\frac{2^3.7}{3.17.103}$  \\\hline \end{tabular}

\begin{tabular}{|c|c|c|c|}\hline $11/2$ & $13/2$ & $15/2$ & $17/2$\\\hline  $\Delta_{25,11}\oplus\Delta_{19}$ & $\Delta_{25,13}^2\oplus\Delta_{19}, \Delta_{25,19,13}$ & $\Delta_{25,15}\oplus\Delta_{19}$ & $\Delta_{25,17}\oplus\Delta_{19}$\\\hline $\frac{2^3}{3^2.5.17}$ & $\frac{2.31}{3^2.5.7.19}$ & $\frac{2^4.3}{5.7^2.19}$ & $\frac{2^5.3}{7.31.43}$\\\hline \end{tabular}
\vskip10pt
Note that this time, for $s=9/2$, there are representations with stable, tempered Arthur parameters, whose existence is not demanded by the appearance of any large prime in the numerator of an $L$-value. Again, the $43$ in the last denominator can be explained as in Example 2. Though this instance of Harder's conjecture appears to have been accidentally omitted from the table at the end of \cite{vdG}, $43$ is a divisor of $L_{\alg}(f,j+k)$, with $f$ of weight $26$. In fact, although this instance of Harder's conjecture has not itself been proved, it happens to be the particular example for which Ibukiyama proved his half-integral weight version in \cite[Theorem 4.4]{I}.

\section{Some more experimental congruences when $n=3$}
So far we have found some congruences, then checked the occurrence of large prime moduli in numerators of ratios of $L$-values. As already noted in Section 5, the congruences should lead to the construction of elements in Selmer groups, then the Bloch-Kato conjecture explains the appearance of the large primes in the ratios of $L$-values. So from a congruence we should predict a factor in an $L$-ratio. But the conjecture in \S 3.3 actually goes in the opposite direction, saying that a large prime should occur in an $L$-value only as a result of a congruence. In the previous section, various large primes showed up other than the ones we were looking for, so a good test of the conjecture would be now to find experimental evidence for the congruences which conjecturally follow from this. We have put this experimental evidence in this later section to emphasize this logical point.

{\bf Example 1: $\mathbf{\Delta_{25,15,9}, q=557}$.} This arose in Example 2 of the previous section. The first $\tr(T(p))$ in the table is on a $3$-dimensional space, but subtracting off an endoscopic contribution gives
$$T(p)(\Delta_{25,15,9})=\tr(T(p))-[2p^5T(p)(\Delta_{15})+\tr(T(p))(\Delta_{25,9}^2)],$$
where $\Delta_{25,9}^2$ corresponds to $(j,k)=(8,10)$, for which the space of level $1$ genus $2$ cusp forms is $2$-dimensional.
\vskip10pt
\begin{tabular}{|c|c|c|c|c|}\hline $p$ & $\tr(T(p))$ & $T(p)(\Delta_{15})$ & $\tr(T(p))(\Delta_{25,9}^2)$ & $T(p)(\Delta_{25,15,9})$ \\\hline $2$ & $15216$ & $216$ & $7440$ & $-6048$\\\hline $3$ & $-557532$ & $-3348$ & $1348560$ & $-278964$ \\\hline $5$ & $717423510$ & $52110$ & $-141412200$ & $533148210$\\\hline $7$ & $64935299016$ & $2822456$ & $-22882568800$ & $-7056168168$\\\hline $11$ & $9763224800748$ & $20586852$ & $448932567408$ & $2683226030436$ \\\hline\end{tabular}

\begin{tabular}{|c|c|c|c|}\hline $p$ & $T(p)(\Delta_{25,15,9})$ & $T(p)(\Delta_{25,15})$ & $-T(p)(\Delta_{25,15,9})$\\ & & & $+[T(p)(\Delta_{25,15})+p^{8}+p^{17}]$\\\hline $2$ & $-6048$ & $-3696$ & $2^4.3.5.\mathbf{557}$\\\hline $3$ & $-278964$ & $511272$ & $2^6.3^6.5.\mathbf{557}$\\\hline $5$ & $533148210$ & $118996620$ & $2^8.3.5.17.67.313.\mathbf{557}$\\\hline $7$ & $-7056168168$ & $-82574511536$ & $2^7.3^3.5.293.\mathbf{557}.82463$\\\hline $11$ & $2683226030436$ & $5064306707064$ & $2^6.3.5.7.17.211.\mathbf{557}.37646261$\\\hline\end{tabular}
\vskip10pt
The data is consistent with the predicted congruence $$T(p)(\Delta_{25,15,9})\equiv T(p)(\Delta_{25,15})+p^{8}+p^{17}\pmod{557},$$
in fact we have checked the congruence for all primes $p\leq 53$.

{\bf Example 2: $\mathbf{\Delta_{25,19,13}, q=31}$.} This arose in Example 5 of the previous section. The first $\tr(T(p))$ in the table is on a $3$-dimensional space, but subtracting off an endoscopic contribution,
$$T(p)(\Delta_{25,19,13})=\tr(T(p))-[2p^3T(p)(\Delta_{19})+\tr(T(p))(\Delta_{25,13}^2)],$$
where $\Delta_{25,13}^2$ corresponds to $(j,k)=(12,8)$, for which the space of level $1$ genus $2$ cusp forms is $2$-dimensional. We have checked the congruence for all primes $p\leq 53$, and show the results for $p\leq 11$.
\vskip10pt
\begin{tabular}{|c|c|c|c|c|}\hline $p$ & $\tr(T(p))$ & $T(p)(\Delta_{19})$ & $\tr(T(p))(\Delta_{25,13}^2)$ & $T(p)(\Delta_{25,19,13})$ \\\hline $2$ & $6432$ & $456$ & $-1536$ & $672$\\\hline $3$ & $2206116$ & $50652$ & $173232$ & $-702324$ \\\hline $5$ & $140035350$ & $-2377410$ & $724983000$ & $9404850$\\\hline $7$ & $2180027592$ & $-16917544$ & $28504729184$ & $-14719266408$\\\hline $11$ & $-1608110653332$ & $-16212108$ & $-24717511671792$ & $23152557649956$\\\hline\end{tabular}

\begin{tabular}{|c|c|c|c|}\hline $p$ & $T(p)(\Delta_{25,19,13})$ & $T(p)(\Delta_{25,19})$ & $-T(p)(\Delta_{25,19,13})$\\ & & & $+[T(p)(\Delta_{25,19})+p^{6}+p^{19}]$\\\hline $2$ & $672$ & $-2880$ & $2^5.3.5^2.7.\mathbf{31}$\\\hline $3$ & $-702324$ & $-538920$ & $2^8.3^3.5^2.7.\mathbf{31}^2$\\\hline $5$ & $9404850$ & $118939500$ & $2^{10}.3.5^2.7,\mathbf{31}^2.36919$\\\hline $7$ & $-14719266408$ & $1043249200$ & $2^9.3^4.5^2.7^3.13.\mathbf{31}.79537$\\\hline $11$ & $23152557649956$ & $-9077287359096$ & $2^8.3.5^2.7.\mathbf{31}.706711927.20771$\\\hline\end{tabular}
\vskip10pt

{\bf Example 3: $\mathbf{\Delta_{25,17,3}^2, q=59}$.} This arose in Example 1 of the previous section. We have a $2$-dimensional $M(V_{\mu},K)$, but we might find the eigenvalues $a,b$ of $T(p)$ by solving the quadratic equation $x^2-(a+b)x+ab=0$, where $a+b=\tr(T(p))$ and $ab=\frac{1}{2}((a+b)^2-(a^2+b^2))=\frac{1}{2}((\tr(T(p)))^2-\tr(T(p)^2))$.
The problem becomes to find $\tr(T(p)^2)$.

In the language of \S 2, $T(p)=T_{f_1}=p^{a_1-(5/2)}T'_{f_1}=p^{10}T'_{f_1}$. Similarly we define $T(p^2):=p^{20}T'_{2f_1}$, $T(p,p):=p^{20}T'_{f_1+f_2}$. In the local Hecke algebra at $p$ there is a relation
$${T'}_{f_1}^2=T'_{2f_1}+(p+1)T'_{f_1+f_2}+(p^5+p^4+p^3+p^2+p+1),$$
i.e. $$T(p)^2=T(p^2)+(p+1)T(p,p)+p^{a_1-5}(p^5+p^4+p^3+p^2+p+1),$$
where here $a_1=25$. Such relations may be proved using various ideas expounded by Gross in \cite{G}. See \cite[\S 7.1.3]{Me} and \cite[VI, Exemple 2.11]{CL} for something similar.
We have the traces of $T(4)$ and $T(2,2)$, as well as of $T(2)$, so we can calculate the trace of $T(2)^2$.
\vskip10pt
\begin{tabular}{|c|c|c|c|c|}\hline $\tr(T(2))$ & $\tr(T(4))$ & $\tr(T(2,2))$ & $\tr(T(2)^2)$ & $\frac{1}{2}((\tr(T(2)))^2-\tr(T(2)^2))$ \\\hline $-768$ & $-36421632$ & $-29859840$ & $6119424$ & $-2764800$\\\hline
\end{tabular}

\begin{tabular}{|c|c|}\hline  $T(2)(\Delta_{25,17,3}^2)$ & $-T(2)(\Delta_{25,17,3}^2)+T(2)(\Delta_{25,17})+2^{11}+2^{14}$\\\hline $-384\pm 192\sqrt{79}$ & $22416\mp 192\sqrt{79}$\\\hline
\end{tabular}
\vskip10pt
For the predicted congruence to be possible, we need this difference to have norm divisible by $59$, and we find indeed that
$\Norm(22416\mp 192\sqrt{79})=2^8.3^3.5^2.7^2.\mathbf{59}$.
\vskip10pt
{\bf Example 4: $\mathbf{\Delta_{25,17,7}^2, q=1223}$.} This arose in Example 1 of the previous section. This time we have a $3$-dimensional $M(V_{\mu},K)$, but after subtracting an endoscopic contribution from $\tr(T(2))$ and $\tr(T(2)^2)$ we may proceed as above.
\vskip10pt
\begin{tabular}{|c|c|c|c|c|}\hline $\tr(T(2))$ & $\tr(T(4))$ & $\tr(T(2,2))$ & $\tr(T(2)^2)$ & $T(2)(\Delta_{25,7}\oplus\Delta_{17})$\\\hline $-14832$ & $79978752$ & $65968128$ & $476064000$ & $-11616+2^4(-528)=-20064$\\\hline
\end{tabular}

\begin{tabular}{|c|c|c|c|}\hline $\tr(T(2)|_{\Delta_{25,17,7}^2})$ & $\tr(T(2)^2|_{\Delta_{25,17,7}^2})$ & $T(2)(\Delta_{25,17,7}^2)$ & $-T(2)(\Delta_{25,17,7}^2)$\\ & & & $+T(2)(\Delta_{25,17})+2^9+2^{16}$\\\hline $5232$ & $73499904$ & $2616\pm 216\sqrt{641}$ & $67032\mp 216\sqrt{641}$\\\hline \end{tabular}
\vskip10pt
$\Norm(67032\mp 216\sqrt{641})=2^{12}.3^4.11.\mathbf{1223}$.
\vskip10pt
{\bf Example 5: $\mathbf{\Delta_{25,19,5}^2, q=103}$.} This arose in Example 5 of the previous section, and is similar to the previous example.
\vskip10pt
\begin{tabular}{|c|c|c|c|c|}\hline $\tr(T(2))$ & $\tr(T(4))$ & $\tr(T(2,2))$ & $\tr(T(2)^2)$ & $T(2)(\Delta_{25,5}\oplus\Delta_{19})$\\\hline $10176$ & $3207168$ & $-22394880$ & $134203392$ & $-96+2^3(456)=3552$\\\hline
\end{tabular}

\begin{tabular}{|c|c|c|c|}\hline $\tr(T(2)|_{\Delta_{25,19,5}^2})$ & $\tr(T(2)^2|_{\Delta_{25,19,5}^2})$ & $T(2)(\Delta_{25,19,5}^2)$ & $-T(2)(\Delta_{25,19,5}^2)$\\ & & & $+T(2)(\Delta_{25,19})+2^9+2^{16}$\\\hline $6624$ & $121586688$ & $3312\pm 240\sqrt{865}$ & $27600\mp 240\sqrt{865}$\\\hline \end{tabular}
\vskip10pt
$\Norm(27600\mp 240\sqrt{865})=2^{11}.3^3.5^3.\mathbf{103}$.

\section{The setup for $G=\SO(n,n), M\simeq \GL_2\times\SO(n-2,n-2)$}
We use similar notation to \S 2. Let
$$G=\SO(n,n)=\{g\in M_{2n}:\,\,^tg{J}g={J},\,\det(g)=1\},$$
where
$${J}=\begin{pmatrix}0_n & I_n\\I_n & 0_n\end{pmatrix}.$$
This is a connected, semi-simple algebraic group, split over  $\Q$.
It has a maximal torus $T= \{\diag(t_1,\ldots,t_n,t_1^{-1},\ldots,t_n^{-1}) : t_1,\ldots,t_n \in \GL_1\}$ with character group $X^*(T)$ spanned by $\{e_1,\ldots,e_n\}$ where $e_i$ sends $\diag(t_1,\ldots,t_n,t_1^{-1},\ldots,t_n^{-1})$ to $t_i$ for $1\leq i\leq n$. The cocharacter group $X_*(T)$ is spanned by $\{f_1,\ldots,f_n\}$, where $f_1:t\mapsto \diag(t,1,\ldots,1,t^{-1},1,\ldots,1)$, etc.~and so $\langle e_i,f_j\rangle=\delta_{ij}$, where $\langle,\rangle:X^*(T)\times X_{*}(T)\rightarrow\Z$ is the natural pairing. We can order the roots so that the set of positive roots is $\Phi_G^+=\{e_i-e_j:\,i<j\}\cup \{e_i+e_j:\,i< j\}$, with simple positive roots $\Delta_G=\{e_1-e_2,e_2-e_3,\ldots,e_{n-1}-e_n,e_{n-1}+e_n\}$. The half-sum of the positive roots is $\rho_G=(n-1)e_1+(n-2)e_2+\cdots +e_{n-1}$. The Weyl group $W_G$ is generated by permutations of the $t_i$ and by inversions swapping $t_i$ with $t_i^{-1}$. The long element $w_0^G$ is the product of all the inversions.

Suppose now that $n\geq 3$. If we choose the simple root $\alpha=e_2-e_3$, this determines a maximal parabolic subgroup $P=MN$, where $N$ is the unipotent radical and $M$ is the Levi subgroup, characterised by $\Delta_M=\Delta_G-\{\alpha\}$, and then $M\simeq \GL_2\times\SO(n-2,n-2)$. The positive roots occurring in the Lie algebra of $N$ are $\Phi_N=\Phi^+_G-\Phi^+_M=\{e_1-e_3,\ldots,e_1-e_n,e_1+e_2,\ldots,e_1+e_n, e_2-e_3,\ldots,e_2-e_n,e_2+e_3,\ldots,e_2+e_n\}$, i.e.~those positive roots whose expression as a sum of simple roots includes $\alpha$. The half-sum is $\rho_P=\frac{2n-3}{2}(e_1+e_2)$, and $\langle\rho_P,\check{\alpha}\rangle=\frac{2n-3}{2}$, where $\check{\alpha}$ is the coroot associated with $\alpha$. Let $\tilde{\alpha}:=\frac{1}{\langle\rho_P,\check{\alpha}\rangle}\rho_P=e_1+e_2$. The Langlands dual group $\hat{G}\simeq G$.

Let $\Pi_f$ be the unitary, cuspidal automorphic representation of $\GL_2(\A)$ associated to a cuspidal Hecke eigenform $f$ of weight $k$ and trivial character, and let $\Pi'$ be a unitary, cuspidal, automorphic representation of $\SO(n-2,n-2)(\A)$. Let $\Pi=\Pi_f\times\Pi'$, which is a unitary, cuspidal, automorphic representation of $M(\A)$. Let $\lambda=\frac{k-1}{2}(e_1-e_2)+a_1e_3+\cdots a_{n-2}e_n$, with $a_1\geq a_2\geq \ldots\geq a_{n-2}\geq 0$, be the infinitesimal character of $\Pi_{\infty}$, up to $W_M$. We shall assume that the $a_i$ are all distinct, from each other and from $(k-1)/2$, with $a_{n-2}>0$.

 For any prime $p$ such that the local component $\Pi_p$ is unramified, let $\chi_p=-[\log_p(\alpha_p)(e_1-e_2)+\log_p(\beta_1)e_3+\log_p(\beta_2)e_3+\cdots\log_p(\beta_{n-2})e_n]\in X^*(T)\otimes_{\Z}\C$ be such that $\Pi_p$ is isomorphic to the (unitarily) parabolically induced representation $\Ind_{B(\QQ_p)}^{M(\QQ_p)}(|\chi_p|_p)$, where $B$ is a Borel subgroup of $M$ containing $T$. Note that $p^{(k-1)/2}(\alpha_p+\alpha_p^{-1})=a_p(f)$, where $f=\sum_{i=1}^{\infty}a_n(f)q^n$, with $a_1(f)=1$. This $\chi_p\in X^*(T)\otimes i\RR$ gives rise to a Satake parameter $t(\chi_p)\in \hat{T}(\CC)\subset \hat{M}(\CC)$, as before.

The adjoint representation $r:\hat{M}\rightarrow \Aut(\hat{\n})$, is $r_1\oplus r_2$, with $\Phi_N^1=\{e_1\pm e_{j+2}, e_2\pm e_{j+2}: 1\leq j\leq n-2\}$ and $\Phi_N^2=\{e_1+e_2\}$.
\vskip10pt
\begin{tabular}{|c|c|c|c|}\hline $\gamma\in\Phi_N$ & $\check{\gamma}$ & $\langle\lambda+s\tilde{\alpha},\check{\gamma}\rangle$ &  $|\chi_p(\check{\gamma}(p))|_p$\\\hline  $e_1-e_{j+2}$ ($1\leq j\leq n-2$) & $f_1-f_{j+2}$ & $\frac{k-1}{2}-a_{j}+s$ & $\alpha_p\beta_{j}^{-1}$\\$e_2-e_{j+2}$ ($1\leq j\leq n-2$) & $f_2-f_{j+2}$ & $-\frac{k-1}{2}-a_{j}+s$ & $\alpha_p^{-1}\beta_{j}$\\$e_1+e_{j+2}$ ($1\leq j\leq n-2$) & $f_1+f_{j+2}$ & $\frac{k-1}{2}+a_{j}+s$ & $\alpha_p\beta_{j}^{-1}$\\$e_2+e_{j+2}$ ($1\leq j\leq n-2$) & $f_2+f_{j+2}$ & $-\frac{k-1}{2}+a_{j}+s$ & $\alpha_p^{-1}\beta_{j}$\\$e_1+e+2$ & $f_1+f_2$ & $2s$ & $1$\\\hline \end{tabular}
\vskip10pt
Using the table, $L_{\Sigma}(s,\Pi,r_1)$ is the $L$-function associated with $\Pi_f\times\Pi'$ and the tensor product of the standard representations of $\GL_2$ and $\SO(n-2,n-2)$, while $L_{\Sigma}(s,\Pi,r_2)=\zeta_{\Sigma}(s)$.

For $s>0$, the representation $\Ind_P^G(\Pi\otimes|s\tilde{\alpha}|)$ of $G(\A)$ has infinitesimal character (at $\infty$) $\lambda+s\tilde{\alpha}$ (up to $W_G$-action). We need $s\in \frac{1}{2}+\ZZ$ for $L_{\Sigma}(1+2s,\Pi,r_2)$ to be critical, then for all the $a_i$ to be in $\frac{1}{2}+\Z$, for $\lambda+s\tilde{\alpha}$ to be algebraically integral.

Let $1\leq t\leq n-2$ be such that $\frac{k-1}{2}$ is in between $a_t$ and $a_{t+1}$ (or $\frac{k-1}{2}>a_1$ if $t=1$, $\frac{k-1}{2}<a_{n-2}$ if $t=n-2$).
$$\lambda+s\tilde{\alpha}=\left(\frac{k-1}{2}+s\right)e_1+\left(-\frac{k-1}{2}+s\right)e_2+a_1e_3+\cdots a_{n-2}e_n.$$ Then, for the obvious choice of $w\in W_G$,
$$w(\lambda+s\tilde{\alpha})=a_1e_1+\cdots +a_te_t$$ $$+\left(\frac{k-1}{2}+s\right)e_{t+1}+\left(\frac{k-1}{2}-s\right)e_{t+2}+a_{t+1}e_{t+3}+\cdots +a_{n-2}e_n,$$ which is dominant and regular if we add the condition $s<\min\{a_t-\frac{k-1}{2},\frac{k-1}{2}-a_{t+1}\}$ to those already imposed. This coincides with the condition for $L_{\Sigma}(1+s,\Pi,r_1)$ to be critical. We exclude the smallest value $s=1/2$ from the conjecture below.

Suppose that $q>2\max\langle\lambda,\check{\gamma}\rangle +1=k+2a_1$, and let $\q\mid q$ be a prime divisor of $\q$ in a number field sufficiently large to accommodate all the Hecke eigenvalues and normalised $L$-values we shall consider.

The main conjecture of \cite{BD} is that if $\ord_{\q}(L_{\alg,\Sigma}(1+is,\Pi,r_i)>0$ then there exists a tempered, cuspidal, automorphic representation $\tilde{\Pi}$ of $G(\A)$, unramified outside $\Sigma$, and with $\tilde{\Pi}_{\infty}$ of infinitesimal character $w(\lambda+s\tilde{\alpha})$, such that for all $p\notin\Sigma$, and all $\mu\in X_*(T)$, the eigenvalues
of $T_{\mu}$ on $\tilde{\Pi}_p$ and $\Ind_P^G(\Pi_p\otimes|s\tilde{\alpha}|_p)$ are congruent modulo $\q$.

The standard representation of $\hat{G}$ has highest weight $f_1$ (identifying $X^*(\hat{T})$ with $X_{*}(T)$) and complete set of weights $\{\pm f_1,\pm f_2,\ldots, \pm f_n\}$. Given that this is a single $W_G$-orbit, i.e. that $f_1$ is a minuscule weight, we can calculate the ``right-hand-side'' of the congruence in the following way.
The Satake parameter of $\Ind_P^G(\Pi_p\otimes|s\tilde{\alpha}|_p)$ is $-[\log_p(\alpha_p)(e_1-e_2)+\log_p(\beta_1)e_3+\log_p(\beta_2)e_3+\cdots\log_p(\beta_{n-2})e_n]+s(e_1+e_2)$. Using this,
\vskip10pt
\begin{tabular}{|c|c|}\hline $\mu$ & $|(\chi_p+s\tilde{\alpha})(\mu(p))|_p$ \\ \hline $f_1$ & $\alpha_pp^{-s}$\\$f_2$ & $\alpha_p^{-1}p^{-s}$\\$f_{i+2}$ ($1\leq i\leq n-2$) & $\beta_i$ \\ \hline \end{tabular}
\vskip10pt
The trace is $(\alpha_p+\alpha_p^{-1})(p^s+p^{-s})+\sum_{i=1}^{n-2}(\beta_i+\beta_i^{-1})$. We multiply by $p^{\langle w(\lambda+s\tilde{\alpha}),f_1\rangle}$ to get the eigenvalue for $T_{f_1}$:
$$T_{f_1}(\Pi_p\otimes |s\tilde{\alpha}|_p)$$ $$=\begin{cases} a_p(f)(1+p^{2s})+\sum_{i=1}^{n-2}p^{(k-1)/2+s}(\beta_i+\beta_i^{-1}) & \text{ if $\frac{k-1}{2}>a_1$};\\(p^{(a_1-(k-1)/2)+s}+p^{(a_1-(k-1)/2)-s})a_p(f)+\sum_{i=1}^{n-2}p^{a_1}(\beta_i+\beta_i^{-1}) & \text{ if $\frac{k-1}{2}<a_1$}.\end{cases}$$

\section{The case $n=4$ with $i=1$}
As above, we have $\Pi_f$ the unitary, cuspidal, automorphic representation of $\GL_2(\A)$ associated to a cuspidal Hecke eigenform $f$ of weight $k$, trivial character, and $\Pi'$ a unitary, cuspidal, automorphic representation of $\SO(2,2)(\A)$. From now on we assume that $\Pi_f$ and $\Pi'$ are unramified at all finite $p$. As in the proof of \cite[Proposition* 4.15]{CR}, using the central isogeny $\SO(2,2)\rightarrow\PGL(2)\times\PGL(2)$, we can get $\Pi'$ by giving a pair $\Pi_g, \Pi_h$ of cuspidal, automorphic representations of $\PGL(2)(\A)$, associated to cuspidal Hecke eigenforms $g$, $h$ of level $1$, let's say of weights $\ell, m$ respectively. (Strictly speaking, $\Pi'$ is a discrete automorphic representation, but by \cite[Theorem 4.3]{Wa} it will be cuspidal.)
The infinitesimal character of $\Pi'$ is $\frac{\ell+m-2}{2}e_3+\frac{|\ell-m|}{2}e_4$, and for each finite prime $p$, its Satake parameter is the tensor product of those of $\Pi_f$ and $\Pi_g$, so that the standard $L$-function of $\Pi'$ is $L(\Pi_g\otimes\Pi_h,s)$, and $L(s,\Pi,r_1)$ is the triple product $L$-function $L(\Pi_f\otimes\Pi_g\otimes\Pi_h,s)$.

We relabel $\{k,\ell,m\}=\{k_1,k_2,k_3\}$ in such a way that $k_1\geq k_2\geq k_3$, and we also relabel $\{f,g,h\}=\{f_{k_1},f_{k_2},f_{k_3}\}$ in the obvious way. Henceforth we consider only examples for which $k_1<k_2+k_3$ and for which each of $f_{k_1}, f_{k_2}$ and $f_{k_3}$ spans its space of cusp forms. The field of coefficients will then be $\Q$. We have
$$L(\Pi_f\otimes\Pi_g\otimes\Pi_h,s)=L(s+\frac{k_1+k_2+k_3-3}{2},f\otimes g\otimes h).$$
If,
$$\hat{L}(s,f\otimes g\otimes h):=\Gamma_{\C}(s)\Gamma_{\C}(s-(k_1-1))\Gamma_{\C}(s-(k_2-1))\Gamma_{\C}(s-(k_3-1))L(s,f\otimes g\otimes h),$$
where $\Gamma_{\C}(s)=2(2\pi)^{-s}\Gamma(s)$, then $\hat{L}(s)=- \hat{L}(k_1+k_2+k_3-2-s)$. Let
$$\hat{L}_{\alg}(s,f\otimes g\otimes h):=\frac{\hat{L}(s,f\otimes g\otimes h)}{(f,f)(g,g)(h,h)},$$
with Petersson norms in the denominators. For integers $k_1\leq t\leq k_2+k_3-2$, $L(t,f\otimes g\otimes h)$ is a critical value. The condition $q>k+2a_1=k_1+k_2+k_3$ appearing in the conjecture on congruences is too large for some of our examples. We can do better by viewing the conjecture as saying that whenever the Bloch-Kato conjecture predicts that a $q$-torsion Selmer group is non-trivial, it gets that way via a mod $q$ congruence of Hecke eigenvalues. The bound on $q$ guarantees that the $q$-parts of Tamagawa factors (in the appropriate normalisation) are trivial, so that when we see a factor $q$ in the numerator of a normalised $L$-value, Bloch-Kato accounts for it in the Selmer group. But for the special form of tensor product motive here, we can prove more about the Tamagawa factors, and thus employ smaller bounds for $q$, as in the proposition below.

As in \cite[Lemma 5.1]{DH}, the product of Petersson norms is the normalised Deligne period, up to a power of $(2\pi i)$ that depends on the point of evaluation and is taken care of by the $\Gamma_{\C}$ factors. Note in \cite[(5.2)]{DH}, which, following Hida, notes the congruence factor $c(f_{k_i})$ intervening between $(f_{k_i},f_{k_i})$ and $\Omega_{f_{k_i}}^+\Omega_{f_{k_i}}^-$, that $c(f_{k_i})$ is trivial for us, because the spaces are $1$-dimensional, leaving no room for the congruences measured by the $c(f_{k_i})$. (Even if it were non-trivial, it would only contribute to the denominator of $\hat{L}_{\alg}(t)$ anyway.) These periods are well-defined up to primes less than $k_1$. Suppose that $q>k_1$, and that $k_1\leq t\leq k_2+k_3-2$ is an integer, avoiding the central point $t=\frac{k_1+k_2+k_3-2}{2}$. Then the Bloch-Kato conjecture predicts that
$$\ord_{q}\left(\hat{L}_{\alg}(t)\right)$$ $$=\ord_{q}\left(\frac{c_q(t)\#H^1_{f}(\Q,T_{q}^*(1-t)\otimes (\Q_{q}/\Z_{q}))}{\#H^0(\Q,T_{q}^*(1-t)\otimes(\Q_{q}/\Z_{q}))\#H^0(\Q,T_{q}(t)\otimes(\Q_{q}/\Z_{q}))}\right),$$
where $c_q(t)$ is a certain Tamagawa factor. We shall not define it here, but it can be dealt with as in \cite[Proposition 7.5]{DH}, following \cite[Proposition 2.16]{DFG}. In \cite{DH} the weights are all equal, making things slightly simpler, but the idea is essentially the same, so we merely state the result.

\begin{prop} Suppose that $f_{k_1}$ is ordinary at $q$, that $k_1\leq t'\leq k_2+k_3-2$ and $q>\max\{k_1,2k_3-2-(t'-(k_1-1)),k_3+2+t'-k_2\}$. Then $\ord_q(c_q(t'))\leq 0$.
\end{prop}
This is most effective for $t'$ left of the central point, whereas we are primarily interested in $t=k_1+k_2+k_3-2-t'=1+s+\frac{k_1+k_2+k_3-3}{2}$, to the right of the central point.
But if $\ord_{q}\left(\hat{L}_{\alg}(t)\right)>0$ then by the functional equation, $\ord_{q}\left(\hat{L}_{\alg}(t')\right)>0$. If $q>\max\{k_1,2k_3-2-(t'-(k_1-1)),k_3+2+t'-k_2\}$ and $f_{k_1}$ is ordinary at $q$, Bloch-Kato predicts that $H^1_{f}(\Q,T_{q}^*(1-t')\otimes (\Q_{q}/\Z_{q}))$ is non-trivial, hence, by \cite{Fl}, that $H^1_{f}(\Q,T_{q}^*(1-t)\otimes (\Q_{q}/\Z_{q}))$ is non-trivial. We then predict a congruence
$$T(p)(\tilde{\Pi})\equiv (p^{((\ell+m-k-1)/2)+s}+p^{((\ell+m-k-1)/2)-s})a_p(f)+a_p(g)a_p(h)\pmod{q},$$
for all primes $p$, where $T(p)=T_{f_1}$ and $\tilde{\Pi}$ is a tempered, cuspidal, automorphic representation of $\SO(4,4)(\A)$, with $\tilde{\Pi}_{\infty}$ of infinitesimal character
$$\frac{\ell+m-2}{2}e_1+\left(\frac{k-1}{2}+s\right)e_2+\left(\frac{k-1}{2}-s\right)e_3+\frac{|\ell-m|}{2}e_4.$$

To compute Hecke eigenvalues for $\tilde{\Pi}$ we proceed as described in \S 6, except now we use $\SO(8)$, the orthogonal group of the $E_8$ lattice, with $\SO(8)(\Z)\simeq W(E_8)^+$, of order $348364800$. For $\mu=a_1e_1+a_2e_2+a_3e_3+a_4e_4$ (in the notation of \cite[5.2]{CR}, in particular $a_1$ now means something different), with $a_1,a_2,a_3,a_4\in\ZZ$ and $a_1\geq a_2\geq a_3\geq a_4\geq 0$, let $V_{\mu}$ be the space of the complex representation $\theta_{\mu}$ of $\SO(8)$ with highest weight $\mu$, and let $\rho:=3e_1+2e_2+e_3$. The infinitesimal character of the representation $\theta_{\mu}$ of $\SO(8)(\RR)$ is $\mu+\rho$. Let $K$ be the open compact subgroup $\prod_p \SO(8)(\ZZ_p)$ of $\SO(8)(\A_f)$, and let
$$M(V_{\mu},K):=\{f:\SO(8)(\A_f)\rightarrow V_{\mu}:$$
$$ ~f(gk)=f(g)~\forall k\in K,~f(\gamma g)=\theta_{\mu}(\gamma)(f(g))~\forall \gamma\in \SO(8)(\QQ)\}$$
be the space of $V_{\mu}$-valued algebraic modular forms with level $K$.
Since $$\#(\SO(8)(\QQ)\backslash \SO(8)(\A_f)/K)=1,$$ $M(V_{\mu},K)$ can be identified with the fixed subspace $V_{\mu}^{\SO(8)(\ZZ)}$.
For each (finite) prime $p$, $\SO(8)(\QQ_p)\simeq\SO(4,4)(\QQ_p)$, and the local Hecke algebras are naturally isomorphic. The third named author has computed the trace of $T(p):=T_{f_1}$ on $M(V_{\mu},K)$ for all $p\leq 23$, and for $a_1\leq 12$ \cite{Me, Me2}. In the examples below, we have verified the expected congruences for all primes $p\leq 23$, but we display the data just as far as $p=13$.

{\bf Example 1: $(k,\ell,m)=(18,12,20)$. }
From the computations of Ibukiyama and Katsurada in the appendix, we see that $\hat{L}_{\alg}(26,f\otimes g\otimes h)=\frac{2^{53}.\mathbf{31}}{3.17}$. (Note that our $\hat{L}_{\alg}$ is their $L_{\alg}$.) With $(k_1,k_2,k_3)=(20,18,12)$, the critical range is $20\leq t\leq 28$, and $26=\frac{k_1+k_2+k_3-3}{2}+1+s$ for $s=3/2$. With $t=26$, $t'=22$, $\max\{k_1,2k_3-2-(t'-(k_1-1)),k_3+2+t'-k_2\}=\max\{20,19,18\}=20$, a bound comfortably exceeded by $q=31$. Since $31\nmid -104626880141728=a_{31}(f_{20})$, $f_{20}$ is ordinary at $31$. According to \cite[Table 9]{CR}, there is a single stable, tempered Arthur parameter $\Delta_{30,20,14,8}$ for the relevant infinitesimal character $15e_1+10e_2+7e_3+4e_4$, and a table at \cite{CRtab} shows that $\dim(V_{\mu}^{\SO(8)(\ZZ)})=1$ (for $\mu=12e_1+8e_2+6e_3+4e_4$), so we obtain $T(p)(\Delta_{30,20,14,8})$ as the trace of $T(p)$ on $V_{\mu}^{\SO(8)(\ZZ)}$.
\vskip10pt
\begin{tabular}{|c|c|c|c|}\hline $p$ & $a_p(f_{18})$ & $a_p(f_{12})$ & $a_p(f_{20})$\\\hline $2$ & $-528$ & $-24$ & $456$\\$3$ & $-4284$ & $252$ & $50652$\\$5$ & $-1025850$ & $4830$ & $-2377410$\\$7$ & $3225992$ & $-16744$ & $-16917544$\\$11$ & $-753618228$ & $534612$ & $-16212108$\\$13$ & $2541064526$ & $-577738$ & $50421615062$\\\hline\end{tabular}

\begin{tabular}{|c|c|c|}\hline $p$ & $T(p)(\Delta_{30,20,14,8})$ & $(p^5+p^8)a_p(f_{18})+a_p(f_{12})a_p(f_{20})-T(p)(\Delta_{30,20,14,8})$\\\hline
$2$ & $69120$ & $-2^6.3^2.13.\mathbf{31}$\\$3$ & $10614240$ & $-2^9.3^5.7.\mathbf{31}$\\$5$ & $-18486732600$ & $-2^{10}.3^2.5^2.7.17.\mathbf{31}.467$\\$7$ & $-984888553600$ & $2^{11}.3^4.7^3.23.\mathbf{31}.491$\\$11$ & $-4326973699452192$ & $-2^9.3^2.5^2.17 .\mathbf{31}.317.8175953$\\$13$ & $-59262235173721720$ & $2^{10}.3^4.7^2.\mathbf{31}.617.27064319$\\\hline\end{tabular}
\vskip10pt

{\bf Example 2: $(k,\ell,m)=(22,12,20)$. }
From the computations of Ibukiyama and Katsurada in the appendix, we see that $\hat{L}_{\alg}(29,f\otimes g\otimes h)=\frac{2^{57}.3.7.\mathbf{73}}{5.19}$. With $(k_1,k_2,k_3)=(22,20,12)$, the critical range is $22\leq t\leq 30$, and $29=\frac{k_1+k_2+k_3-3}{2}+1+s$ for $s=5/2$. With $t=29$, $t'=23$, $\max\{k_1,2k_3-2-(t'-(k_1-1)),k_3+2+t'-k_2\}=\max\{22,20,17\}=22$, a bound comfortably exceeded by $q=73$. Since $73\nmid -43284759511102937494=a_{73}(f_{22})$, $f_{22}$ is ordinary at $73$. By \cite[Table 9]{CR}, there is a single stable, tempered Arthur parameter $\Delta_{30,26,16,8}$ for the relevant infinitesimal character $15e_1+13e_2+8e_3+4e_4$, and a table at \cite{CRtab} shows that $\dim(V_{\mu}^{\SO(8)(\ZZ)})=1$ (for $\mu=12e_1+11e_2+7e_3+4e_4$), so we obtain $T(p)(\Delta_{30,26,16,8})$ as the trace of $T(p)$ on $V_{\mu}^{\SO(8)(\ZZ)}$.
\vskip10pt
\begin{tabular}{|c|c|c|}\hline $p$ & $a_p(f_{22})$ & $T(p)(\Delta_{30,26,16,8})$\\\hline
$2$ & $-288$ & $-80496$\\$3$ & $-128844$ & $12133152$\\$5$ & $21640950$ & $-28999867896$\\$7$ & $-768078808$ & $6809124360320$\\$11$ & $-94724929188$ & $-2979055414026720$\\$13$ & $-80621789794$ & $-66466630034660152$\\\hline\end{tabular}

\begin{tabular}{|c|c|}\hline $p$ & $(p^2+p^7)a_p(f_{22})+a_p(f_{12})a_p(f_{20})-T(p)(\Delta_{30,26,16,8})$\\\hline $2$ & $2^4.3^3.\mathbf{73}$\\$3$ & $-2^7.3^4.\mathbf{73}.373$\\$5$ & $2^8.3^3.\mathbf{73}.109.31069$\\$7$ & $-2^9.3^4.5.\mathbf{73}.211.401.499$\\$11$ & $-2^7.3^3.\mathbf{73}.277.26371814347$\\$13$ & $-2^8.3^4.\mathbf{73}.3317356371521$\\\hline
\end{tabular}
\vskip10pt

{\bf Example 3: $(k,\ell,m)=(22,12,18)$. }
From the computations of Ibukiyama and Katsurada in the appendix, we see that $\hat{L}_{\alg}(27,f\otimes g\otimes h)=\frac{2^{52}.3.\mathbf{43}}{19}$. With $(k_1,k_2,k_3)=(22,18,12)$, the critical range is $22\leq t\leq 28$, and $27=\frac{k_1+k_2+k_3-3}{2}+1+s$ for $s=3/2$. With $t=27$, $t'=23$, $\max\{k_1,2k_3-2-(t'-(k_1-1)),k_3+2+t'-k_2\}=\max\{22,20,19\}=22$, a bound comfortably exceeded by $q=43$. Since $43\nmid -193605854685795844=a_{43}(f_{22})$, $f_{22}$ is ordinary at $43$. According to \cite[Table 9]{CR}, there is a single stable, tempered Arthur parameter $\Delta_{28,24,18,6}$ for the relevant infinitesimal character $14e_1+12e_2+9e_3+3e_4$, and a table at \cite{CRtab} shows that $\dim(V_{\mu}^{\SO(8)(\ZZ)})=1$ (for $\mu=11e_1+10e_2+8e_3+3e_4$), so we obtain $T(p)(\Delta_{28,24,18,6})$ as the trace of $T(p)$ on $V_{\mu}^{\SO(8)(\ZZ)}$.
\vskip10pt
\begin{tabular}{|c|c|c|}\hline $p$ & $T(p)(\Delta_{28,24,18,6})$ & $(p^2+p^5)a_p(f_{22})+a_p(f_{12})a_p(f_{18})-T(p)(\Delta_{28,24,18,6})$\\\hline $2$ & $-10080$ & $2^5.3^2.\mathbf{43}$\\$3$ & $3900960$ & $-2^9.3^5.7.\mathbf{43}$\\$5$ & $1700332200$ & $2^{10}.3^2.5^2.7.\mathbf{43}.887$\\$7$ & $-95141488000$ & $-2^{11}.3^3.\mathbf{43}.991.5477$\\$11$ & $-50025639432672$ & $-2^9.3^2.5^2.\mathbf{43}.877.3595481$\\$13$ & $-1259590157649880$ & $-2^{10}.3^3.7.\mathbf{43}^2.3209.26261$\\\hline
\end{tabular}
\vskip10pt

{\bf Example 4: $(k,\ell,m)=(16,12,20)$. }
From the computations of Ibukiyama and Katsurada in the appendix, we see that $\hat{L}_{\alg}(26,f\otimes g\otimes h)=\frac{2^{54}.3.5.\mathbf{19}}{13.17}$. With $(k_1,k_2,k_3)=(20,16,12)$, the critical range is $20\leq t\leq 26$, and $26=\frac{k_1+k_2+k_3-3}{2}+1+s$ for $s=5/2$. With $t=26$, $t'=20$, $\max\{k_1,2k_3-2-(t'-(k_1-1)),k_3+2+t'-k_2\}=\max\{20,21,18\}=21$, a bound not quite achieved by $q=19$, but we shall look for the congruence anyway. Note that even in $\hat{L}_{\alg}(20)$, the factor $19$ only appears in the numerator because of the factor $\Gamma(20)=19!$, i.e. it does not appear in the numerator of $L_{\alg}(20)$. But by the same token it appears in the denominators of the other critical values $L_{\alg}(t)$ for $21\leq t\leq 25$, suggesting that it appears in the denominators of the corresponding Tamagawa factors $c_{19}(t)$. Therefore it would not be too surprising if it occurred also in the denominator of $c_{19}(26)$, allowing $\#H^1_{f}(\Q,T_{q}^*(1-t))$ still to be non-trivial.

 According to \cite[Table 9]{CR}, there is a single stable, tempered Arthur parameter $\Delta_{30,20,10,8}$ for the relevant infinitesimal character $15e_1+10e_2+5e_3+4e_4$, and a table at \cite{CRtab} shows that $\dim(V_{\mu}^{\SO(8)(\ZZ)})=1$ (for $\mu=12e_1+8e_2+4e_3+4e_4$), so we obtain $T(p)(\Delta_{30,20,10,8})$ as the trace of $T(p)$ on $V_{\mu}^{\SO(8)(\ZZ)}$.
\vskip10pt
\begin{tabular}{|c|c|c|}\hline $p$ & $a_p(f_{16})$ & $T(p)(\Delta_{30,20,10,8})$\\\hline
$2$ & $216$ & $52992$\\$3$ & $-3348$ & $7306848$\\$5$ & $52110$ & $671424840$\\$7$ & $2822456$ & $-107393799808$\\$11$ & $20586852$ & $-167258251753632$\\$13$ & $-190073338$ & $-42627620077539832$\\\hline\end{tabular}

\begin{tabular}{|c|c|}\hline $p$ & $(p^5+p^{10})a_p(f_{16})+a_p(f_{12})a_p(f_{20})-T(p)(\Delta_{30,20,10,8})$\\\hline $2$ & $2^6.3^3.5.\mathbf{19}$\\$3$ & $-2^9.3^4.5.7^2.\mathbf{19}$\\$5$ & $2^{10}.3^3.5.\mathbf{19}.37.5113$\\$7$ & $2^{11}.3^4.5.\mathbf{19}.193.262271$\\$11$ & $2^9.3^3.5.13.\mathbf{19}.8550379.3659$\\$13$ & $-2^{10}.3^4.5.\mathbf{19}.863.271279.14197$\\\hline
\end{tabular}
\vskip10pt

{\bf Example 5: $(k,\ell,m)=(20,12,16)$. }
We are looking at the same $L$-value as in the previous example, but $f,g$ and $h$ have been permuted. Again, $19$ is not strictly speaking big enough, but we look for the congruence anyway. By \cite[Table 9]{CR}, there is a single stable, tempered Arthur parameter $\Delta_{26,24,14,4}$ for the relevant infinitesimal character $13e_1+12e_2+7e_3+2e_4$, and a table at \cite{CRtab} shows that $\dim(V_{\mu}^{\SO(8)(\ZZ)})=1$ (for $\mu=10e_1+10e_2+6e_3+2e_4$), so we obtain $T(p)(\Delta_{26,24,14,4})$ as the trace of $T(p)$ on $V_{\mu}^{\SO(8)(\ZZ)}$.
\vskip10pt
\begin{tabular}{|c|c|c|}\hline $p$ & $T(p)(\Delta_{26,24,14,4})$ & $(p+p^6)a_p(f_{20})+a_p(f_{12})a_p(f_{16})-T(p)(\Delta_{26,24,14,4})$\\\hline $2$ & $-16128$ & $2^4.3^3.5.\mathbf{19}$\\$3$ & $-1851552$ & $2^9.3^3.5.\mathbf{19}.29$\\$5$ & $313754760$ & $-2^{10}.3^3.5.\mathbf{19}.37.383$\\$7$ & $34598801792$ & $-2^{11}.3^4.5.\mathbf{19}.131497$\\$11$ & $-25141764069792$ & $2^9.3^3.5.\mathbf{19}.5655173$\\$13$ & $232075615185608$ & $2^{10}.3^4.5^2.\mathbf{19}.643.1613.5953$\\\hline
\end{tabular}
\vskip10pt
In the appendix there are several further examples of primes $q>k_1$ dividing numerators of normalised $L$-values, but not only do these examples exceed the $a_1\leq 12$ for which traces have been computed, they also have $\dim(V_{\mu}^{\SO(8)(\ZZ)})>1$. For example, with $(k,\ell,m)=(20,16,18)$ we have
$\hat{L}_{\alg}(30,f\otimes g\otimes h)=\frac{2^{58}.7.\mathbf{2297}}{13.17}$. The relevant $\mu=13e_1+11e_2+5e_3+e_4$, and $\dim(V_{\mu}^{\SO(8)(\ZZ)})=4$. For another example, with $(k,\ell,m)=(22,16,20)$ we have
$\hat{L}_{\alg}(31,f\otimes g\otimes h)=\frac{2^{59}.3.7.\mathbf{6619}}{13.17.19}$. The relevant $\mu=14e_1+11e_2+7e_3+2e_4$, and $\dim(V_{\mu}^{\SO(8)(\ZZ)})=12$. So we have not attempted to test the predicted congruences.
\section{The special case $f=h$.}

Since $k_1<k_2+k_3$, we have $k-1>\frac{\ell -1}{2}$.
Define integers $a,b,c$ by $a+3=k-1, b+2=\frac{\ell-1}{2}+s$ and $c+1=\frac{\ell-1}{2}-s$. For $\frac{1}{2}<s<\min\{\frac{\ell+m-2}{2}-\frac{k-1}{2},\frac{k-1}{2}-\frac{|\ell-m|}{2}\}$, considering separately the cases $k\leq \ell$ and $k\geq \ell$ (with $k=m$), one checks that $\frac{\ell-1}{2}+s<k-1$, so that $a\geq b\geq c\geq 0$. Note that $k=m=a+4$, $\ell=b+c+4$ and $t=a+b+6$, with $s=\frac{b-c+1}{2}$ and $\frac{k+\ell+m-3}{2}=\frac{2a+b+c+9}{2}$.

Seventeen examples of $\hat{L}_{\alg}(a+b+6,f\otimes f\otimes g)$ appear in \cite[Table 3]{IKPY}. Actually, since some of their examples involve weights other than $12,16,18,20,22$ and $26$, it is actually the norm of this algebraic number that appears in their table. Their computations
are connected with the seventeen experimental congruences supporting \cite[Conjecture 10.8]{BFvdG}, which is also discussed in \cite[Section 8, Case 2]{BD}, where, as in \cite{IKPY}, $f$ and $g$ are the other way round. The right-hand-sides of these congruences are $a_p(f)(a_p(g)+p^{b+2}+p^{c+1})$, which is exactly what we get in the case $f=h$. The left-hand-sides are Hecke eigenvalues for genus $3$ vector-valued Siegel cusp forms, of ``type'' $(a,b,c)$. These can be equated with the desired Hecke eigenvalues of cuspidal, automorphic representations of $\SO(4,4)(\A)$, of infinitesimal character $\frac{a+b+c+6}{2}e_1+\frac{a+b-c+4}{2}e_2+\frac{a-b+c+2}{2}e_3+\frac{|a-(b+c)|}{2}e_4$, via the conjectured functorial lift from $\PGSp_3$ to $\SO(4,4)$, associated with the homomorphism $\Spin(4,3)\rightarrow \SO(4,4)$ of $L$-groups that is the $8$-dimensional spinor representation. Thus, the congruences of the previous section may be viewed as a generalisation of \cite[Conjecture 10.8]{BFvdG} (which is due to the authors of that paper in collaboration with Harder and Mellit).

Actually, $L(a+b+6,f\otimes f\otimes g)=L(a+b+6,\Sym^2 f\otimes g)L(b+3,g)$. In the seventeen examples referred to above, the modulus of the congruence always appears in the first factor (suitably normalised), in fact it is only that factor that actually appears in \cite[Conjecture 10.8]{BFvdG}. What about the other factor? If $\q\mid q$ with $q>b+c+4$, and $\ord_{\q}(L_{\alg}(b+3,g))>0$, then according to Harder's conjecture there should exist a cuspidal Hecke eigenform $F$ of genus $2$, level $1$, vector-valued of type $\Sym^{b-c}\otimes\det^{c+3}$, with Hecke eigenvalue at $p$ congruent to $a_p(g)+p^{b+2}+p^{c+1}$ mod $\q$. Let $\Pi_F$ be the associated cuspidal, automorphic representation of $\PGSp_2(\A)$. Then the conjectured congruence (with right hand side $a_p(f)(a_p(g)+p^{b+2}+p^{c+1})$) would be satisfied by the conjectured functorial lift of $\Pi_f\times\Pi_F$ from $\PGL_2\times\PGSp_2$ to $\SO(4,4)$, via the homomorphism of $L$-groups
$$\SL(2)\times\Spin(3,2)\rightarrow\SL(2)\times\Sp_2\rightarrow\SO(4,4),$$ where $\Spin(3,2)\rightarrow \Sp_2$ is the spinor representation and the second arrow is the ``tensor product'' representation.

Recalling \S 4, we should also consider the case $i=2$, no longer assuming that $f=h$. This may be done using an endoscopic lift from
$\SO(2,2)(\A)\times\SO(2,2)(\A)$ to $\SO(4,4)(\A)$. One representation of $\SO(2,2)(\A)$ comes via tensor product lift from $g$ and $h$, the other similarly from $f$ and something satisfying a Ramanujan-style congruence.

\appendix
\section{Some triple $L$ values \\ by Tomoyoshi Ibukiyama and Hidenori Katsurada}
The aim of this short note is to give some tables of rigorous explicit critical values
of the triple $L$ functions.
A program to give explicitly any critical values of the triple $L$ functions
rigorously, or norms of the values when the value is not rational,
was prepared by authors in 2011 during writing up a paper \cite{IKPY}
on triple $L$ values and related congruences between modular forms of higher
degree jointly written with C. Poor and D. Yuen.
 There are already many explicit examples in
 \cite{IKPY} including the norms of the
 values when they are not rational, but
 this note is an expanded version. Theoretically nothing is
 new here and we run the same program, but the paper \cite{IKPY}
 treated various other aspects and has a complicated appearance,
 so it seems not useless to
 give a short explanation focussed only on the critical values.
 Although there exists another program in Magma
to calculate these kinds of values by numerical
approximation based on Dokchitser \cite{Do},
our program is based on completely different theory and
gives rigorous values for the  algebraic parts of any critical values
of triple $L$ functions, i.e. the critical values
divided by the Petersson inner products of the three forms (and elementary
Gamma factors). For a theoretical explanation, see \cite{boechererschulzepillot} and
\cite{IKPY}, which will
be outlined also below.

The (right part of the) critical points of the triple $L$ function of
elliptic modular forms $f,g,h$ of weights $k_1$, $k_2$, $k_3$ are given by
\[
\frac{k_1+k_2+k_3+r}{2}-2
\]
where $r$ is a positive integer with $r\geq 2$, also bounded above as described later.
(The functional equation is $s\rightarrow k_1+k_2+k_3-2-s$.)
In actual calculation we use the Siegel Eisenstein series $E_r$ of degree $3$ of even
weight $r$ and its pullback formula which is explained later.
For degree $3$, the original Siegel Eisenstein series itself
converges only for weight $r>4$, but
we can define holomorphic Eisenstein series $E_{r}$ for even weight $r\geq 2$
by Hecke's trick by analytic continuation of the real analytic Eisenstein
series.
Here, when $r=2$, then the analytically continued Eisenstein series
vanishes identically (For these facts,  see \cite{shimura}, \cite{haruki} for example).
On the other hand, because the sign in the functional equation is minus, the $L$-function, which is analytic at the central point
$s=(k_1+k_2+k_3)/2-1$, must vanish there. (See \cite{tsatoh}, \cite{mizumoto},
 \cite{boechererschulzepillot}).
So we may assume that $r$ is an even integer with
$r\geq 4$ and consider the critical values only for these.

We also assume the three weights of elliptic modular forms which define
the triple $L$ function to be in the so-called balanced case. That is, we assume that
there exist non-negative integers $\nu_1$, $\nu_2$, $\nu_3$ and $r\geq 2$
such that
\begin{align}
k_1 & = r+\nu_2+\nu_3, \notag \\
k_2 & = r+\nu_3+\nu_1, \quad \label{weight}\\
k_3 & = r+\nu_1+\nu_2. \notag
\end{align}
So we have
\begin{align*}
\nu_1 & =(k_2+k_3-k_1-r)/2, \\
\nu_2 & = (k_3+k_1-k_2-r)/2, \\
\nu_3 & = (k_1+k_2-k_3-r)/2.
\end{align*}
The non-negativity of $\nu_i$ means that we are assuming that
$k_1<k_2+k_3$ if $k_1\geq k_2 \geq k_3$.
Also this gives the upper bound $r\leq k_2+k_3-k_1$,
which matches the bound of the critical points.
If we take $k_1\geq k_2\geq k_3$, then $\nu_1\leq \nu_2\leq \nu_3$.
These numbers determine the necessary differential operators
which are used in the pullback formula of Eisenstein series of degree 3.
In proving this formula we followed B\"ocherer and Schulze-Pillot, who looked at $r=2$ (the central point) and level $N>1$ in \cite[Theorem 5.7]{boechererschulzepillot}, but also implicitly treated the other $r$ in \cite[(2.1),(2.41)]{boechererschulzepillot}. By means of \cite[Proposition 4.5]{IKPY}, we replaced their differential operators by those from \cite{ibuzagier}, which are more amenable to computation.

Roughly speaking, the pullback formula in our case is explained as follows.
We have a differential operator
${\mathbb D_{r,\nu_1.\nu_2,\nu_3}}$ which depends on $r$ and $\nu_i$, and which acts on
any Siegel modular form $F$ of weight $r$
on the Siegel upper half space $H_3$ of degree $3$.
The restriction of ${\mathbb D}_{r,\nu_1,\nu_2,\nu_3}F$ to the diagonal elements
$H_1\times H_1 \times H_1$ of $H_3$ gives the elliptic modular forms of
weight $k_1$, $k_2$, $k_3$ respectively for the diagonal variables,
where $k_1$, $k_2$, $k_3$ are determined by
the above described relations. If we take $F=E_{r}$, then the diagonal restriction of
${\mathbb D}_{r,\nu_1,\nu_2,\nu_3}E_{r}$ is a linear combination of
tensors of elliptic eigenforms of weights $k_1$, $k_2$ and $k_3$. The critical $L$-value $L\left(\frac{k_1+k_2+k_3+r}{2}-2,f\otimes g\otimes h\right)$, up to elementary factors and Petersson norms, appears in the coefficient of $f(z_1)g(z_2)h(z_3)$. See \cite[Theorem 4.8]{IKPY}. So the only problem is to write down
the diagonal restriction of ${\mathbb D}_{r,\nu_1,\nu_2,\nu_3}E_{r}$
as a linear combination of such products.
We need three things in order to calculate the explicit critical values.
We fix weights $k_1$, $k_2$, $k_3$ and consider all integers
$r\geq 4$, $\nu_i\geq 0$ such that (\ref{weight}) holds. Note again that the choice of $r$, hence $E_r$, determines which critical value we are looking at. Then we need \\
(i) explicit Fourier coefficients of the Eisenstein series of degree three
of weight $r$; \\
(ii) explicit differential operators ${\mathbb D}_{r,\nu_1,\nu_2,\nu_3}$; \\
(iii) explicit elliptic modular forms. \\
For (i), a complete closed formula is known in \cite{katsurada}. For (ii),
a complete explicit formula for such differential operators, taken from \cite{ibuzagier}, immediately precedes \cite[Theorem 4.3]{IKPY}, and (iii) is classical.

As examples, we write down below all the critical values for
three forms of different weights belonging to
$S_{k}(\SL(2,\Z))$ in the case that these spaces are one dimensional.
For calculation, this condition is not necessary, but for the sake of simplicity,
we assumed this here.
This means that $k=12$, $16$, $18$, $20$, $22$, or $26$,
and here we assume that the three forms are different combinations of these.
Other explicit examples of values,
in cases where two weights are equal, sometimes norms of the values
for forms in spaces of dimension greater than one, have been given
in \cite{IKPY}. There is no difference in the method of calculation.
Note that here, since $k_1>k_2>k_3$, necessarily $r<k_3$, so only cusp forms appear on the right-hand-side of the pullback formula \cite[Theorem 4.8]{IKPY}, and since additionally the spaces are $1$-dimensional, there is only a single term.

Now since the normalization of the $L$ functions might sometimes differ
in different contexts, we give here the precise definition of the algebraic part of the
$L$ function we adopted.
For elliptic modular forms $f$, $g$, $h$ of weight $k_1$, $k_2$, $k_3$
with $k_3\leq k_2\leq k_1<k_2+k_3$ and
for integers $l$ with
\[
\frac{k_1+k_2+k_3}{2}-1 \leq l \leq k_2+k_3-2 \qquad
\]
we define as in \cite{IKPY} the algebraic part $L_{\alg}(l)$ of the triple $L$ values
$L(l)$ for $f$, $g$, $h$ by
\begin{multline*}
L_{\alg}(l,f\otimes g \otimes h)
\\ = \frac{L(l,f\otimes g \otimes h)\Gamma_{\C}(l)\Gamma_{\C}(l-k_1+1)
\Gamma_{\C}(l-k_2+1)\Gamma_{\C}(l-k_3+1)}{(f,f)(g,g)(h,h)},
\end{multline*}
where $\Gamma_{\C}(s)=2(2\pi)^{-s}\Gamma(s)$ and
$(f,f)$, $(g,g)$ $(h,h)$ are Petersson inner metric of the forms defined
by the metric $y^{k-2}dxdy$ on the fundamental domain.
In other words, we have
\[
\dfrac{L(l,f\otimes g\otimes h)}
{\pi^{4l+3-k_1-k_2-k_3}(f,f)(g,g)(h,h)}
=
\dfrac{2^{4l-1-k_1-k_2-k_3}L_{\alg}(l,f\otimes g\otimes h)}
{(l-1)!(l-k_1)!(l-k_2)!(l-k_3)!}.
\]

We examine an example of the ratio for different $l$ in order to
check that  it is consistent with the calculation of Magma.
We have
\begin{align*}
\frac{L(l_1,f\otimes g \otimes h)}{\pi^{4(l_1-l_2)}L(l_2,f\otimes g \otimes h)}
& =\frac{L_{\alg}(l_1,f\otimes g \otimes h)}{L_{\alg}(l_2,f\otimes g \otimes h)}
 \\& \times
\frac{2^{4(l_1-l_2)}
\Gamma(l_2)\Gamma(l_2-k_1+1)\Gamma(l_2-k_2+1)\Gamma(l_2-k_3+1)}
{\Gamma(l_1)\Gamma(l_1-k_1+1)\Gamma(l_1-k_2+1)\Gamma(l_1-k_3+1)}
\end{align*}
For example, if $(k_1,k_2,k_3)=(20,16,12)$ and $l_1=26$ and $l_2=24$, then we have
\begin{align*}
\frac{L(26,f\otimes g \otimes h )}{\pi^{8}L(24,f\otimes g \otimes h)}
& =
\frac{2^{8}\Gamma(24)\Gamma(5)\Gamma(9)\Gamma(13)
L_{\alg}(26,f \otimes g \otimes h)}
{\Gamma(26)(\Gamma(7)\Gamma(11)\Gamma(15)L_{\alg}(24,f \otimes g \otimes h)}
\\ & =\frac{2^2\cdot L_{\alg}(26,f \otimes g \otimes h)}
{3^4\cdot 5^4 \cdot 7 \cdot 13 \cdot L_{\alg}(24,f \otimes g \otimes h)}.
\end{align*}
By our program, we have
\begin{align*}
L_{\alg}(24,f\otimes g \otimes h) & = -2^{52}\cdot 7 \cdot 13^{-1} \cdot 17^{-1}, \\
L_{\alg}(25,f\otimes g \otimes h) & = -2^{53}\cdot 3 \cdot 17^{-1}, \\
L_{\alg}(26,f\otimes g \otimes h) & =  -2^{54}\cdot 3 \cdot 5 \cdot 13^{-1},
 \cdot 17^{-1} \cdot 19,
\end{align*}
so for example we have
\[
\frac{L(26,f\otimes g \otimes h)}{\pi^8 L(24,f\otimes g \otimes h)}
=
\frac{2^4\cdot 19}{3^3\cdot 5^3 \cdot 7^2 \cdot 13}.
\]
This coincides with the value produced by Magma.

We denote by $f_{k}$ the normalized
cusp form belonging to $S_{k}(\SL(2,\Z))$ such that $\dim S_{k}(\SL(2,\Z))=1$.
Then we have $f_{12}  = q-24 q^2 + \cdots$, $f_{16}  = q+216 q^2 + \cdots$,
$ f_{18}  = q-528 q^2+\cdots $, $f_{20}  = q+456q^2 + \cdots $,
$f_{22}  = q-288q^2+\cdots $ and $f_{26} = q-48q^2+\cdots$. These
coefficients are needed for our calculation. \\

Now we give tables of critical values below.

{\setstretch{1.2}
\[
\begin{array}{|c|c|c|}\hline
r & \text{critical points $s$} & L_{\alg}(s,f_{12}\otimes f_{16}\otimes f_{18})
\\ \hline
4 & 23 &  2^{45}\cdot 3\cdot 5 / 13 \\ \hline
6 & 24 &  2^{48}\cdot 3\cdot 5 / 13 \\ \hline
8 & 25 &  2^{47}\cdot 3^3 \\ \hline
10 & 26 & 2^{52}\cdot 3^3\cdot 7 / 13  \\ \hline
\end{array}
\]
}

{\setstretch{1.2}
\[
\begin{array}{|c|c|c|}\hline
r & \text{critical points $s$} & L_{\alg}(s,f_{12}\otimes f_{16}\otimes f_{20})
\\ \hline
4 & 24 & 2^{52}\cdot 7 / (13\cdot 17) \\ \hline
6 & 25 & 2^{53} \cdot 3 / 17 \\ \hline
8 & 26 & 2^{54}\cdot 3 \cdot 5 \cdot 19/ (13 \cdot 17) \\ \hline
\end{array}
\]
}

{\setstretch{1.2}
\[
\begin{array}{|c|c|c|}\hline
r & \text{critical points $s$} & L_{\alg}(s,f_{12}\otimes f_{16}\otimes f_{22})
\\ \hline
4 & 25 & 2^{51}\cdot 3^2 \cdot  7 / (5 \cdot 19)  \\ \hline
6 & 26 & 2^{54}\cdot 3^3\cdot  7 / (13 \cdot19)  \\ \hline
\end{array}
\]
}

{\setstretch{1.2}
\[
\begin{array}{|c|c|c|}\hline
r & \text{critical points $s$} & L_{\alg}(s,f_{12}\otimes f_{16}\otimes f_{26})
\\ \hline
2 & 26 & 0  \\ \hline
\end{array}
\]
}

{\setstretch{1.2}
\[
\begin{array}{|c|c|c|}\hline
r & \text{critical points $s$} & L_{\alg}(s,f_{12}\otimes f_{18}\otimes f_{20})
\\ \hline
4 &  25 & 2^{50}\cdot 5 \cdot  11 / (7 \cdot 17)  \\ \hline
6 &  26 & 2^{53}\cdot  31/(3\cdot 17)  \\ \hline
8 &  27 & 2^{52}\cdot 3^2\cdot 5\cdot 7 /17 \\ \hline
10 & 28  & 2^{59}\cdot 3\cdot 5^3/ (7\cdot 17) \\ \hline
\end{array}
\]
}

{\setstretch{1.2}
\[
\begin{array}{|c|c|c|}\hline
r & \text{critical points $s$} & L_{\alg}(s,f_{12}\otimes f_{18}\otimes f_{22})
\\ \hline
4 & 26  &  2^{52} \cdot 11 / 19 \\ \hline
6 & 27 &  2^{52}\cdot 3 \cdot 43/19 \\ \hline
8 & 28 &  2^{56}\cdot 3^3 \cdot 5 /19 \\ \hline
\end{array}
\]
}

{\setstretch{1.2}
\[
\begin{array}{|c|c|c|}\hline
r & \text{critical points $s$} & L_{\alg}(s,f_{12}\otimes f_{18}\otimes f_{26})
\\ \hline
4 & 28 &  2^{61}\cdot 3 \cdot 5 /(7\cdot 23) \\ \hline
\end{array}
\]
}

{\setstretch{1.2}
\[
\begin{array}{|c|c|c|}\hline
r & \text{critical points $s$} & L_{\alg}(s,f_{12}\otimes f_{20}\otimes f_{22})
\\ \hline
4 &  27 &  2^{55}\cdot 11 \cdot 13 / (3\cdot 17 \cdot 19) \\ \hline
6 &  28 &  2^{58}\cdot 11^2 / (17 \cdot 19) \\ \hline
8 &  29 &  2^{57}\cdot 3 \cdot 7 \cdot 73 / (5 \cdot 19) \\ \hline
10 & 30 & 2^{62}\cdot 3^3\cdot 7^2 /(5\cdot 19) \\ \hline
\end{array}
\]
}

{\setstretch{1.2}
\[
\begin{array}{|c|c|c|}\hline
r & \text{critical points $s$} & L_{\alg}(s,f_{12}\otimes f_{20}\otimes f_{26})
\\ \hline
4 &  29 &  2^{59}\cdot 13 / 23 \\ \hline
6 &  30 &  2^{62}\cdot 3^3 /23 \\ \hline
\end{array}
\]
}

{\setstretch{1.2}
\[
\begin{array}{|c|c|c|}\hline
r & \text{critical points $s$} & L_{\alg}(s,f_{12}\otimes f_{22}\otimes f_{26})
\\ \hline
4 & 30 & 2^{65}\cdot 3\cdot 13/(5\cdot 19 \cdot 23) \\ \hline
6 & 31 & 2^{63}\cdot 3\cdot 13/23 \\ \hline
8 & 32 & 2^{66}\cdot 3^{2}\cdot 5\cdot 31 /(11\cdot 23) \\ \hline
\end{array}
\]
}

{\setstretch{1.2}
\[
\begin{array}{|c|c|c|}\hline
r & \text{critical points $s$} & L_{\alg}(s,f_{16}\otimes f_{18}\otimes f_{20})
\\ \hline
4 & 27 & 2^{53}\cdot 5 \cdot 7/(13\cdot 17) \\ \hline
6 & 28 & 2^{60}\cdot 5/(13\cdot 17) \\ \hline
8 & 29 & 2^{54}\cdot 11 \cdot 719/ (13 \cdot 17) \\ \hline
10 & 30 & 2^{58}\cdot 7 \cdot 2297/ (13\cdot 17) \\ \hline
12 & 31 & 2^{56}\cdot 3 \cdot 5^{2}\cdot 11^{2}\cdot 23/17 \\ \hline
14 & 32 & 2^{60}\cdot 3^3\cdot 5^3\cdot 11\cdot 19/17  \\ \hline
\end{array}
\]
}

{\setstretch{1.2}
\[
\begin{array}{|c|c|c|}\hline
r & \text{critical points $s$} & L_{\alg}(s,f_{16}\otimes f_{18}\otimes f_{22})
\\ \hline
4 & 28 & 2^{57}\cdot 3\cdot 5/(13\cdot 19) \\ \hline
6 & 29 & 2^{55}\cdot 3\cdot 7 \cdot 53/(13\cdot 19) \\ \hline
8 & 30 & 2^{57}\cdot 3^4\cdot 7 \cdot 61 /(5\cdot 13 \cdot 19) \\ \hline
10 &  31  & 2^{57}\cdot 3^2\cdot 7 \cdot 283/19\\ \hline
12 &  32  & 2^{59}\cdot 3^3\cdot 5 \cdot 7 \cdot 11 \\ \hline
 \end{array}
\]
}

{\setstretch{1.2}
\[
\begin{array}{|c|c|c|}\hline
r & \text{critical points $s$} & L_{\alg}(s,f_{16}\otimes f_{18}\otimes f_{26})
\\ \hline
4 &  30 &  2^{60}\cdot 3\cdot 31 /(13\cdot 23) \\ \hline
6 &  31 &  2^{60}\cdot 3^3\cdot 5/23\\ \hline
8 &  32 &  2^{64}\cdot 3^{2}\cdot 5^{2}/23 \\ \hline
 \end{array}
\]
}

{\setstretch{1.2}
\[
\begin{array}{|c|c|c|}\hline
r & \text{critical points $s$} & L_{\alg}(s,f_{16}\otimes f_{20}\otimes f_{22})
\\ \hline
4 & 29 & 2^{61}\cdot 7 \cdot 11/(3\cdot 5\cdot 17 \cdot 19) \\ \hline
6 & 30 &  2^{61}\cdot 7 \cdot 541/(3\cdot 13 \cdot 17 \cdot 19) \\ \hline
8 & 31 &  2^{59}\cdot 3\cdot 7\cdot 6619 /(13\cdot 17 \cdot 19) \\ \hline
10 & 32 & 2^{63}\cdot 7\cdot 31 \cdot 1511/(13 \cdot 17 \cdot 19) \\ \hline
12 & 33 & 2^{61}\cdot 3 \cdot 7^2\cdot 11^2 \cdot 83 / (5\cdot 19) \\ \hline
14 & 34 & 2^{65}\cdot 3^4\cdot 7^2 \cdot 11^3 /(5 \cdot 17) \\ \hline
\end{array}
\]
}

{\setstretch{1.2}
\[
\begin{array}{|c|c|c|}\hline
r & \text{critical points $s$} & L_{\alg}(s,f_{16}\otimes f_{20}\otimes f_{26})
\\ \hline
4 &  31 & 2^{62}\cdot 7 \cdot 11/(17\cdot 23)\\ \hline
6 &  32 &  2^{66}\cdot 3^{3}\cdot 5\cdot 11/(13 \cdot 17 \cdot 23)\\ \hline
8 &  33 & 2^{63}\cdot 18911/(11\cdot 23) \\ \hline
10 &  34 & 2^{67}\cdot 3^{3}\cdot 7^{2}\cdot 59 / (17\cdot 23)  \\ \hline
\end{array}
\]
}

{\setstretch{1.2}
\[
\begin{array}{|c|c|c|}\hline
r & \text{critical points $s$} & L_{\alg}(s,f_{16}\otimes f_{22}\otimes f_{26})
\\ \hline
4 &  32 & 2^{69}\cdot 7 /(5\cdot 19\cdot 23) \\ \hline
6 &  33 &  2^{65}\cdot 7\cdot 113 /(19\cdot 23) \\ \hline
8 &  34 & 2^{68}\cdot 3 \cdot 173 \cdot 479/(5\cdot 13 \cdot 19 \cdot 23)\\ \hline
10 & 35 & 2^{68}\cdot 3^3\cdot 1831/(5\cdot 23) \\ \hline
12 & 36 & 2^{70}\cdot 3^4\cdot 5^2\cdot 7\cdot 11/23 \\ \hline \end{array}
\]
}

{\setstretch{1.2}
\[
\begin{array}{|c|c|c|}\hline
r & \text{critical points $s$} & L_{\alg}(s,f_{18}\otimes f_{20}\otimes f_{22})
\\ \hline
4 &  30 & 2^{63}\cdot 3/ (17\cdot 19) \\ \hline
6 &  31 & 2^{59}\cdot 5^2\cdot 41 /(17\cdot 19) \\ \hline
8 &  32 & 2^{60}\cdot 5\cdot 89\cdot 109/(3\cdot 17 \cdot 19) \\ \hline
10 & 33 &  2^{61}\cdot 5 \cdot 7 \cdot 11^{2}\cdot 13/(3 \cdot 17)\\ \hline
12 & 34 & 2^{62}\cdot 11 \cdot 1237 \cdot 3617/(5\cdot 17 \cdot 19) \\ \hline
14 & 35 & 2^{63}\cdot 3\cdot 11\cdot 13^2\cdot 53 \cdot 353 /(5\cdot 19)\\ \hline
16 & 36 & 2^{68}\cdot 3^{2}\cdot 5\cdot 7^{2}\cdot 11^{2}\cdot 13\cdot 17/19\\ \hline
 \end{array}
\]
}

{\setstretch{1.2}
\[
\begin{array}{|c|c|c|}\hline
r & \text{critical points $s$} & L_{\alg}(s,f_{18}\otimes f_{20}\otimes f_{26})
\\ \hline
4 &  32 & 2^{63}\cdot 5\cdot 97/(7\cdot 17\cdot 23) \\ \hline
6 &  33 & 2^{64}\cdot 5\cdot 601 /(3\cdot 17 \cdot 23) \\ \hline
8 &  34 & 2^{65}\cdot 3\cdot 7 \cdot 907/(17\cdot 23) \\ \hline
10 & 35 & 2^{65}\cdot 3\cdot 13\cdot  10061 /(7\cdot 23)\\ \hline
12 & 36 & 2^{68}\cdot 3^3\cdot 5^{2}\cdot 11\cdot 61/23 \\ \hline
\end{array}
\]
}

{\setstretch{1.3}
\[
\begin{array}{|c|c|c|}\hline
r & \text{critical points $s$} & L_{\alg}(s,f_{18}\otimes f_{22}\otimes f_{26})
\\ \hline
4 & 33 & 2^{65}\cdot 5\cdot 13 /(19 \cdot 23) \\ \hline
6 & 34 & 2^{68}\cdot 3^5/(19 \cdot 23) \\ \hline
8 & 35 & 2^{66}\cdot 3 \cdot 7 \cdot 9839 /(5\cdot 19 \cdot 23) \\ \hline
10 & 36 & 2^{68}\cdot 3^2 \cdot 29 \cdot 97/19 \\ \hline
12 & 37 & 2^{68}\cdot 3^2\cdot 5 \cdot 13\cdot 17 \cdot 223/23 \\ \hline
14 &  38 &  2^{74}\cdot 3^2\cdot 5^3 \cdot 7\cdot  7621/(19 \cdot 23) \\ \hline
\end{array}
\]
}

{\setstretch{1.3}
\[
\begin{array}{|c|c|c|}\hline
r & \text{critical points $s$} & L_{\alg}(s,f_{20}\otimes f_{22}\otimes f_{26})
\\ \hline
4 &  34 & 2^{70}\cdot 1091/(3 \cdot 5 \cdot 17 \cdot 19 \cdot 23) \\ \hline
6 &  35 &  2^{68}\cdot 3 \cdot 193/(17 \cdot 19) \\ \hline
8 &  36 & 2^{71} \cdot 7 \cdot 37511 /(3\cdot 17 \cdot 19 \cdot 23) \\ \hline
10 &  37 & 2^{70}\cdot 3\cdot 7^{4}\cdot 37 \cdot 43 /(17 \cdot 19 \cdot 23)\\ \hline
12 &  38 & 2^{74}\cdot 3 \cdot 98161517/(5\cdot 17 \cdot 19 \cdot 23) \\ \hline
14 &  39 & 2^{71}\cdot 13^2\cdot 884069/23 \\ \hline
16 &  40 & 2^{73}\cdot 3^{4}\cdot 5^{3}\cdot 7 \cdot 13\cdot 5113/23 \\ \hline
\end{array}
\]
}

\end{document}